\documentclass[final,3p,times]{elsarticle}

\usepackage{graphicx}

\usepackage{amssymb}

\usepackage[modulo]{lineno}

\usepackage{amsmath,amsthm}
\usepackage{xcolor}
\usepackage[colorlinks=true, allcolors=blue]{hyperref}
\usepackage{multirow}
\usepackage{array}
\usepackage{caption}
\usepackage{tikz}
\usepackage{booktabs}
\usepackage{algorithm}
\usepackage{algorithmicx}
\usepackage{algpseudocode}
\usepackage{todonotes}
\usepackage{newtxtext,newtxmath}
\usepackage{arydshln}

\theoremstyle{remark} 
\newtheorem{remark}{Remark}

\DeclareSymbolFont{letters}{OML}{cmm}{m}{it}
\DeclareMathAlphabet{\mathcal}{OMS}{cmsy}{m}{n}

\newcommand{\sub}[1]{_\mathit{#1}}

\biboptions{sort&compress}

\journal{JCP}

\begin{document}

\begin{frontmatter}

\title{Multi-Stage Preconditioners for Thermal--Compositional--Reactive Flow in Porous Media}

\author[label1,label2]{Matthias A. Cremon\corref{cor1}}
\address[label1]{Department of Energy Resources Engineering, Stanford University}

\cortext[cor1]{Corresponding author.}

\ead{mcremon@stanford.edu}

\author[label2]{Nicola Castelletto}
\address[label2]{Atmospheric, Earth and Energy Division, Lawrence Livermore National Laboratory}
\ead{castelletto1@llnl.gov}

\author[label2]{Joshua A. White}
\ead{white230@llnl.gov}

\begin{abstract}
We present a family of multi-stage preconditioners for coupled thermal-compositional-reactive reservoir simulation problems. The most common preconditioner used in industrial practice, the Constrained Pressure Residual (CPR) method, was designed for isothermal models and does not offer a specific strategy for the energy equation. For thermal simulations, inadequate treatment of the temperature unknown can cause severe convergence degradation. When strong thermal diffusion is present, the energy equation exhibits significant elliptic behavior that cannot be accurately corrected by CPR's second stage. In this work, we use Schur-complement decompositions to extract a temperature subsystem and apply an Algebraic MultiGrid (AMG) approximation as an additional preconditioning stage to improve the treatment of the energy equation. We present results for several two-dimensional hot air injection problems using an extra heavy oil, including challenging reactive In-Situ Combustion (ISC) cases. We show improved performance and robustness across different thermal regimes, from advection dominated (high P\'eclet number) to diffusion dominated (low P\'eclet number). The number of linear iterations is reduced by 40-85\% compared to standard CPR for both homogeneous and heterogeneous media, and the new methods exhibit almost no sensitivity to the thermal regime.
\end{abstract}


\end{frontmatter}


\section{Introduction}
\label{sec::intro}

Thermal simulation of flow in porous media is used in many applications where the coupling between temperature and the other variables is important, from heavy oil recovery to CO$_2$ sequestration. To access huge reserves of heavy oil \citep{Briggs88} there has been a strong interest in thermal recovery methods such as steam injection, Steam Assisted Gravity Drainage (SAGD) and In-Situ Combustion (ISC) \citep{Prats82,Burger77,Lake89}. These methods provide heat to the reservoir in the form of hot fluids (steam) or generate heat in-situ. ISC processes inject air into a heated reservoir, reaching spontaneous ignition and using the exothermic oxidation reactions to lower the oil viscosity, allowing it to be displaced by the injected fluids.
Although these recovery processes have been successfully applied at many fields, the complexity of the physical phenomena makes predictive numerical simulation extremely challenging.   This work tackles one critical need in particular, the necessity of having a robust and scalable linear solver strategy.

The first thermal-compositional-reactive simulation models were developed in the late 1970s and early 1980s \citep{Burger76,Crookston79,Coats80b,Youngren80,Young83,Rubin85}. Due to limitations on the available computing power at the time, those models typically used a small number of components and simplified physical models for phase behavior, as well as basic solvers. The development of fast and accurate algorithms for flash calculations \citep{Michelsen82a,Michelsen82b,Whitson89} paved the way for integrated simulators capable of simulating the tightly coupled system of non-linear, time dependent Partial Differential Equations (PDEs) \citep{stars,eclipse,Cao02,Lapene10b}. Discretizing these PDEs typically leads to a series of large, non-symmetric, ill-conditioned linear systems. Direct linear solvers \citep{Li05,Kourounis18} can be used, but their memory requirements severely limit the number of components and grid size that can be handled. The industry-standard approach for reservoir simulation is now to use an iterative Krylov subspace solver \citep{Saad03}, such as Generalized Minimum RESidual (GMRES) \citep{Saad86}, paired with a suitable preconditioning scheme. 

Recently, the reservoir simulation community has focused heavily on novel preconditioners for multi-physics problems, particularly coupled flow and geomechanics \citep{White11,Haga12,Gries14,White16,Gaspar17,White19}. Although the physical processes are different, the idea of using approximate Schur-complements to extract easy-to-address subsystems is at the core of all multi-stage preconditioners. In the field of reservoir simulation, most isothermal models use the two-stage Constrained Pressure Residual (CPR) preconditioner introduced in \citet{Wallis83,Wallis85}. Initially designed for Black-Oil models, CPR separates the elliptic and hyperbolic components of the problem and addresses them in separate stages.  We remark that it is common in the reservoir simulation community to refer to pressure as the ``elliptic variable'' and saturation as the ``hyperbolic variable.''  This is a slight abuse of terminology, as it is the underlying differential equations that exhibit elliptic or hyperbolic behavior with respect to these variables, but it does provide a convenient shorthand.  The most common CPR implementation, as described in \citet{Lacroix03} or \citet{Cao05}, uses Algebraic Multigrid (AMG) \citep{Ruge87,Stueben01} for the first stage pressure system and an Incomplete LU (ILU) factorization for the second stage, showing excellent performance also for isothermal compositional models \citep{ThcJia09,Voskov12b,Zhou13,Cusini15}. AMG is well-suited to the elliptic pressure system, while ILU readily corrects local errors associated with the saturation field. 

For thermal compositional simulations, heat conduction through the rock can be a dominant mechanism, leading to CPR convergence degradation. Treating the temperature variable in the ILU second stage poorly approximates long-range temperature coupling. \citet{Li15,Li17} illustrated severe convergence issues using CPR on steam injection and SAGD cases. Their solution, called Enhanced CPR (ECPR), is to retain additional variables in the first stage according to a given heuristic derived from the coupling strength between variables. While an improvement over traditional CPR, this approach has two drawbacks: (1) the expanded first stage subsystem loses any physical meaning, and (2) there is no guarantee the resulting matrix will be well suited for AMG. In another relevant work, \citet{Roy19a,Roy19b} recently presented a detailed study of preconditioning strategies for thermal water injection, including a new Schur-complement expression for the temperature subsystem. They present numerical results from two-phase dead-oil cases with heaters or hot fluid injection and show that an adequate treatment of the temperature improves the robustness, performance and scalability of GMRES solvers.

In this paper we propose a family of multi-stage preconditioners suitable for general thermal reservoir simulation, showing reduced sensitivity to the thermal regime and good convergence behavior. We test these preconditioners on two dimensional, thermal-compositional cases with no reactions, as well as on reactive (combustion) cases. We compare these results to traditional CPR as a baseline. We use thermal-compositional simulation results obtained with the Automatic-Differentiation General Purpose Research Simulator (AD-GPRS \citep{Voskov12}), a reservoir simulation framework for coupled thermal-compositional-mechanics processes \citep{Garipov18}. Another key contribution of this work is a specific ordering strategy for the unknowns that guarantees well-posed primary and secondary systems of equations.

The paper is structured as follows. Section \ref{sec::govEq} presents the governing formulation.  Section \ref{sec:ou} presents an unknown ordering and static-condensation strategy. Section \ref{sec::msp} provides a description of traditional CPR and proposes a new family of CPR-like preconditioners. Section \ref{sec::results} presents numerical results for non-reactive and reactive test cases across different thermal regimes for both homogeneous and heterogeneous media. We conclude in Section \ref{sec::concl} with some avenues for future work.

\section{Problem Statement}
\label{sec::govEq}

For the physical model considered here, the coupled mass and energy conservation equations can be cast as a coupled system of Advection-Diffusion-Reaction (ADR) equations,
\begin{equation} \label{ADR}
 \dfrac{\partial \theta}{\partial t} + \nabla\cdot(\textbf{v}\theta) - \nabla \cdot (D\nabla \theta) - R = 0\,,
\end{equation}
with $\theta$ the conserved variable, $D$ the diffusivity, $\textbf{v}$ the velocity field, and $R$ the source term representing both chemical reactions and well sources/sinks.

\subsection{Mass Conservation Equations}

Mass conservation for each fluid component $c$ across phases $p$, with $n_c$ mobile components and $n_p$ phases, reads
\begin{equation} \label{eq::mass}
 \dfrac{\partial}{\partial t}\left(\phi\sum\limits_{p=1}^{n_p}x_{cp}\rho_pS_p\right) + \nabla \cdot \left(\sum\limits_{p=1}^{n_p}x_{cp}\rho_p\textbf{u}_p\right) + \sum\limits_{p=1}^{n_p} x_{cp}\rho_pq_p + q_c^{r} = 0, \hspace{0.6cm} c = 1,\dots,n_c \,.
\end{equation}
Here, $\phi$ is the porosity; $x_{cp}$ is the mole fraction of component $c$ in phase $p$; $\rho_p$, $S_p$, $\textbf{u}_p$, and $q_p$ are the molar density, saturation, velocity, and volumetric flow rate of phase $p$; and $q_c^{r}$ is the source term from reactions. The equation is dimensional and has units of [mol/day]. We do not consider mass diffusion. In this work, we can also have solid, immobile components in the solid $(s)$ phase, obeying the conservation equation
\begin{equation} \label{eq::mass_solid}
 \dfrac{\partial}{\partial t}\left(\phi c_s\right) + q_s^{r} = 0, \hspace{0.6cm} s = 1,\dots,n_s \,,
\end{equation}
with $c_s$ the mole concentration of solid components $s$ and $q_s^r$ the source term from reactions. 

For mobile species, we compute the fluid phase velocities using the extension of Darcy's law \citep{Darcy56} introduced in \citet{Muskat36}
\begin{equation}  \label{Darcy}
\textbf{u}_b = -\dfrac{k_{rb}}{\mu_b}\textbf{k}\left(\nabla P - \rho_b\textbf{g} \right),
\end{equation}
where $b$ is the phase index, $\textbf{k}$ is the permeability tensor, $P$ is the pressure, $k_{rb}$ is the relative permeability of phase $b$, $\mu_b$ is the viscosity of phase $b$, $\textbf{g}$ is the gravity vector and $\rho_b$ the density of phase $b$.

\subsection{Energy Conservation Equation}

For thermal simulations, we also consider the energy conservation,
\begin{equation} \label{eq::energy}
\dfrac{\partial}{\partial t}\left(\phi\sum\limits_{p=1}^{n_p}U_{p}\rho_pS_p+(1-\phi) \tilde U_r\right)+ \nabla \cdot \left(\sum\limits_{p=1}^{n_p}H_{p}\rho_p\textbf{u}_p\right) - \nabla \cdot \left(\kappa\nabla T\right) + \sum\limits_{p=1}^{n_p} H_p\rho_pq_p + q^{h,r} = 0\,.
\end{equation}
Here, $T$ is the temperature; $\tilde U_r$ is the rock volumetric internal energy; $\kappa$ is the thermal conductivity; $U_{p}$ and $H_{p}$ are the internal energy and enthalpy of phase $p$ \citep{Incropera07}; and $q^{h,r}$ is the source term from reactions. The equation is dimensional and has units of [J/day]. Thermal diffusion (also called conduction) usually cannot be neglected due to the large heat conductivity of the rock matrix, leading to one more term in the energy equation than in the mass equations. We consider the thermal conductivity $\kappa$ a single, constant value corresponding to the rock matrix, due to its much larger contribution compared to the fluid phases thermal conductivity \citep{Prats82}.

\subsection{Local Constraints}

Equations \eqref{eq::mass}, \eqref{eq::mass_solid} and \eqref{eq::energy} form a set of $n_c+n_s+1$ coupled equations. To close the system, we have a set of local constraints. For each mobile component $c$, the fugacity $f$ in all present phases must be equal,
\begin{equation} \label{eq::fug}
f_{cj} = f_{ck},\hspace{0.6cm}\forall j\neq k,\ \ c = 1,\dots,n_c \,,
\end{equation}
where $j$ and $k$ are phases indices, and $f_{cp}$ is the fugacity of component $c$ in phase $p$. Fugacity has units of pressure [bars]. These fugacity constraints give $n_c(n_p-1)$ equations.
In each phase, the sum of all molar fractions must be equal to 1, giving $n_p$ equations,
\begin{equation} \label{eq::phase}
\sum\limits_{c\,=\,1}^{n_c} x_{cp} = 1,\hspace{0.6cm}p = 1,\dots,n_p \,.
\end{equation}
Therefore, equations \eqref{eq::mass}, \eqref{eq::mass_solid}, \eqref{eq::energy}, \eqref{eq::fug} and \eqref{eq::phase} form a complete set of $(n_c+1)n_p+n_s+1$ equations. For convenience, in the rest of the paper we will switch the component indices from $c$ to $i$ and introduce the following notation,
\begin{align*}
x_i & = x_{co} \,, \\
y_i & = x_{cg}  \,,\\
w_i & = x_{cw} \,,
\end{align*} 
indicating the mole fraction of component $i$ in different fluid phases (oil, water, gas).

\begin{remark}
We use the standard notation for water \emph{phase} mole fractions ($w_i$), but note that subscript $_w$ denotes the water \emph{component} index. Therefore $w_w$ is the water component mole fraction in the water phase.
\end{remark}

\subsection{Treatment of Water \& Reactions}

We use a Free-Water flash for the phase behavior calculations. A full review can be found in \citet{Iranshahr09} or \citet{Lapene10}. The assumptions of the Free-Water model imply
\begin{align}  
z_i & = x_iO + y_iV, \hspace{1cm} i = 1,\dots,n_c,\ i \neq w  \,, \label{eq::FW1} \\
z_w & = y_wV + W  \,, \label{eq::FW2}
\end{align}
with $O$, $V$ and $W$ the oil, vapor and water phase molar fractions and $z_i$ the overall mole fraction of component $i$. These equations state that the hydrocarbon components cannot be dissolved in the water phase (Eq. \ref{eq::FW1}) and that the water component cannot be dissolved in the oil phase (Eq. \ref{eq::FW2}). This removes $n_c$ fugacity equality constraints by construction, as well as the water phase constraint equation $\sum w_i = 1$. We remove the $n_c$ unknowns that are identically zero (all $w_i$ except for water, and $x_w$), as well as $w_w$ since it is identically one. The final size of the global system is $n_c(n_p-1)+n_s+1$.

The mass and energy source terms for reactions in Eqs. \eqref{eq::mass}, \eqref{eq::mass_solid} and \eqref{eq::energy} use a standard Arrhenius model \citep{Arrhenius89} and are respectively given by
\begin{align} \label{eq::Arrhenius}
q_c^{r} & = \sum\limits_{r\,=\,1}^{n_r} \left(s^P_{r,c} - s^R_{r,c}\right)A_r \exp\left(\dfrac{-E_{a,r}}{RT}\right)\left[A\right]^{n_{A,r}}P_{\text{O}_2}^{n_{\text{O}_2,r}}, \qquad c = 1,\dots,n_c, \\
q_s^{r} & = \sum\limits_{r\,=\,1}^{n_r} \left(s^P_{r,s} - s^R_{r,s}\right)A_r \exp\left(\dfrac{-E_{a,r}}{RT}\right)\left[A\right]^{n_{A,r}}P_{\text{O}_2}^{n_{\text{O}_2,r}}, \qquad s = 1,\dots,n_s, \\
\label{eq::Arrhenius_h}
q^{h,r} & = \sum\limits_{r\,=\,1}^{n_r} H_rA_r \exp\left(\dfrac{-E_{a,r}}{RT}\right)\left[A\right]^{n_{A,r}}P_{\text{O}_2}^{n_{\text{O}_2,r}},
\end{align}
where $n_r$ the number of reactions, $r$ is the reaction index, $A_r$ is the pre-exponential factor in reaction $r$, $E_{a,r}$ is the activation energy in reaction $r$, $R$ is the ideal gas constant, $n_{A,r}$ is the order of reaction for component $A$ in reaction $r$, $[A]$ is the concentration of the reactant, $n_{\text{O}_2,r}$ is the order of reaction for oxygen in reaction $r$, $P_{\text{O}_2}$ is the partial pressure of oxygen, $s^P_{r,c}$ $\big($or $s^P_{r,s}\big)$ and $s^R_{r,c}$ $\big(\textrm{or } s^R_{r,s}\big)$ are the stoichiometry coefficient for component $c$ $(\textrm{or }s)$ in reaction $r$ as a product (superscript $P$) and a reactant (superscript $R$), and $H_r$ the enthalpy of reaction in reaction $r$.

\subsection{Discretization}

We use the Automatic-Differentiation General Purpose Research Simulator (AD-GPRS \citep{Voskov12}) to solve the system of equations using a finite volume, fully implicit discretization and the natural variables formulation. More details about the reservoir simulation implementation can be found in \citet{Aziz79,Coats80b,Cao02,Voskov09}. We omit details here as they are not the central focus of this work.  At a given time step, we simply note that this discretization leads to a non-linear set of residual equations,
\begin{equation}
r(x) = 0\,,
\end{equation}
for the algebraic vector of grid-cell unknowns $x$.  This system is solved using Newton's method with an Appleyard chopping algorithm \citep{eclipse} to improve convergence robustness. We use limits on the relative variable changes, as well as a check to ensure all the bounded variables remain in the physical range. Given a solution estimate $x^k$, an improved estimate $x^{k+1}$ is determined by
\begin{equation}
\begin{aligned}
&\text{solving}    &&J \delta = - r (x^k) \,, \notag \\
&\text{correcting} &&\delta := \texttt{appleyard}(\delta) \,, \notag \\
&\text{updating}   &&x^{k+1} = x^{k} + \delta \,.
\end{aligned}
\end{equation}
Here, $J = \partial r / \partial x$ is the Jacobian system evaluated at $x^k$ and the vector $\delta$ is the Newton update for time step $k+1$.

\section{Static-Condensation Strategy}
\label{sec:ou}

To lower the computational cost of solving the linear system, we first note that it can be partitioned into \emph{primary} and \emph{secondary} unknowns, as
\begin{equation}
\label{eqn:global}
\begin{bmatrix}
J_{11} & J_{12} \\ J_{21} & J_{22} 
\end{bmatrix}
\begin{bmatrix}
\delta_{1} \\ \delta_{2}
\end{bmatrix} = - 
\begin{bmatrix}
r_{1} \\ r_{2}
\end{bmatrix} \,.
\end{equation}
For a suitable choice of partitioning (described below) the secondary unknowns are only coupled within their own grid cell, but are not coupled across cell boundaries.  The block $J_{22}$ is then block-diagonal, and the secondary unknowns can be readily eliminated through static-condensation.  This leads to a reduced system
\begin{equation}
\label{eqn:reduced}
A \delta_1 = b \,,
\end{equation}
where $A = \left(J_{11} - J_{12}J_{22}^{-1} J_{21}\right)$ is the Schur-complement, and $b = \left(-r_1 + J_{12} J_{22}^{-1} r_2\right)$ is a modified right-hand side.  Note that the system \eqref{eqn:reduced} may be directly assembled, without explicitly forming the global system \eqref{eqn:global}. 

There are several subtleties that must be addressed to create a successful condensation strategy. In particular, phases and components can disappear or reappear in a grid cell as the simulation evolves. If a phase appears or disappears, we expand or shrink the set of variables according to the natural-variable formulation (see \citet{Cao02}). In all cases, the number of primary variables should be $n_c + n_s + 1$, corresponding to the number of globally-coupled equations. There is no unique guidance, however, on the specific primary/secondary partitioning adopted.

It is common practice in compositional simulation to use the ``trace components'' assumption, stating that all components are always present in each cell at least in trace quantities \citep{Voskov12b,Zaydullin17}. Although valid in many conditions, including some thermal processes, that assumption cannot be used as is for cases with reactions: it could lead to infinite reactions if the mass/moles of a reactant cannot truly go to zero. If the trace threshold is included in the reaction terms and boundary conditions, we could get a consistent formulation. However, it would still introduce conditioning problems due to the large number of near-zero terms in the Jacobian.
Moreover, at high temperatures, components are being completely vaporized, which is a key part of the displacement process. Removing the trace assumption will introduce many zero terms in the Jacobian matrix when components are not present. In turn, this makes the selection of primary and secondary variables much more challenging, since rows of zeros can be encountered in blocks that should be invertible. In this section, we present a new, general static-condensation strategy in the context of disappearing components. Although we illustrate the process using our test cases from section \ref{sec::results}, in principle the strategy is applicable to an arbitrary number and type of components.

For a general three-phase flow simulation (ignoring the immobile solid phase) we can have seven different phase states in a cell. Due to the specific nature of the problem studied here, however, we only encounter three cases:

\begin{enumerate}
    \item Gas (G) cell. The water has been vaporized and the oil was either vaporized or burnt. There are then $n_c+n_s+2$ unknowns: pressure, temperature, solid concentrations and vapor mole fractions.
    \item Oil-Gas (OG) cell. The water component has been vaporized and is only present in the gas phase, but the oil phase is present. There are $2n_c+n_s+2$ unknowns, adding the oil mole fractions for all non-water components and the gas saturation.
    \item Oil-Water-Gas (OWG) cell. All mobile phases are present. There are $2n_c+n_s+3$ unknowns, adding the oil saturation.
\end{enumerate}
In all cases, we need to pick an appropriate subset of primary unknowns. A mandatory requirement is that the resulting block $J_{22}$ is invertible. This is achieved by aligning equations with properly chosen unknowns within each grid cell. We also note that it is convenient to avoid pivoting to the extent possible when applying Gaussian-Elimination to these small diagonal blocks.  We therefore try to avoid zero diagonal entries appearing in the alignment strategy.

A few choices are immediately clear.  The solid problem has $n_s$ equations and $n_s$ unknowns.  The simplified model adopted here has only one solid concentration.  We always align the solid conservation equation with the solid concentration, guaranteeing a non-zero diagonal value. Due to the different nature (and scaling) of the energy conservation equation, we always align it with temperature.

We now work on the flow unknowns to satisfy the invertibility requirement. When a component is not present, its conservation equation will only show non-zero values in the columns corresponding to that component's mole fractions ($y_i$ and potentially $x_i$) and those values are 1. The same is true for the fugacity constraints. These considerations lead us to align equations with their corresponding unknowns as much as possible. However, recall that we need to pick $n_c$ primary unknowns in the flow problem. We always want to retain pressure as a primary unknown. The saturation columns in the Jacobian only have non-zero values in the conservation equations, since fugacities and phase constraints do not depend on saturations. If we consider saturations as secondary variables, we can introduce a full row of zeros in $J_{22}$, making it singular. Therefore, if saturations are part of the global set of unknowns (i.e. we have two or three phases present), they need to be selected in the primary unknowns. Conversely, when water is present as a phase, we have a fugacity constraint for water. Since we do not retain the water mole fractions in the free-water context, the water K-value is only a function of $P$, $T$ (both primary variables) and $y_w$. We need to have $y_w$ as a secondary variable to avoid another row of zeros in $J_{22}$ for the 3-phase cells.

With these general considerations in mind, we now describe a recommended partitioning strategy for each of the three possible phase states. For reference, the list of the components used in section \ref{sec::results} is given in Table \ref{tab::comps}.

\begin{table}[htb]
    \centering
    \caption{List of components used in section \ref{sec::results}, ordered according to our static-condensation scheme. The first three indices are of most concern, since they may be aligned with pressure and saturations. Note: C$_{50+f}$ is the portion of C$_{50+}$ that is allowed to react, but the components are identical for thermodynamic properties.}
    {\footnotesize
    \begin{tabular}{ccccccccccccc}
    \toprule
        Index & 1 & 2 & 3 & 4 & 5 & 6 & 7 & 8 & 9 & 10 & 11 & 12  \\
    \midrule
       Comp. & N$_2$ & C$_{50+}$ & H$_2$O & C$_{17}$-C$_{21}$ & C$_{22}$-C$_{27}$ & C$_{28}$-C$_{35}$ & C$_{36}$-C$_{49}$ & C$_2$-C$_{11}$  & C$_{50+f}$ & CO$_2$ & O$_2$ & C$_{12}$-C$_{16}$ \\
    \bottomrule
    \end{tabular}
    }
    \label{tab::comps}
\end{table}

\subsection{Gas Cell}

The pure gas case happens upstream of the combustion front, where the water has been vaporized and moved downstream, and the oil has been burnt. In our implementation, we then only have one secondary equation, the vapor phase constraint $\sum y_i = 1$. There are no saturations in the set of flow unknowns, but we retain pressure. We can pick any of the vapor mole fraction as the sole secondary unknown, since it will be a non-zero value in the phase constraint equation.  We do, however, need to be careful that this component does not disappear, since it would lead to a zero value on the diagonal for the first row (given the lack of dependence on pressure). We align pressure with the nitrogen conservation equation since it is the least likely component to disappear in a pure gas cell.

\subsection{Oil-Gas Cell}

Oil-Gas (OG) cells are located between the combustion front and the water front. The temperature is high enough that water can only be present in the vapor phase, or can be absent altogether. If it is absent, we need to align it with its own mole fraction $y_w$, so we identify water as the third component ($y_3 = y_w$). As all of the oil has not yet been burned or displaced, the flow unknowns also include the gas saturation $S_g$ and the oil mole fractions, $x = [x_1,x_2,x_4,\dots,x_{n_c}]$. Note that we do not have $x_3 = x_w$ since the water component cannot be present in the oil phase. We pick pressure, gas saturation, the vapor mole fraction of water and $n_c-3$ oil mole fractions as primary unknowns, and all other $n_c-1$ vapor mole fractions plus the remaining two oil mole fractions as secondary unknowns. We choose $[x_4,\dots,x_{n_c}]$ for the primary block, so that we can align them with their corresponding components. In an OG cell, the least likely component to disappear is the heaviest hydrocarbon component (that cannot react, in our case C$_{50+}$). We align it with the gas saturation to make sure that we do not encounter a zero-diagonal value. We align the secondary equations with their respective fugacity equations, and the two remaining mole fractions (one vapor and one oil) to their respective phase constraints.

\subsection{Oil-Water-Gas Cell}

Finally, Oil-Water-Gas cells are downstream of the water condensation front, in the cold zone. All three phases are present, and we now add the oil saturation $S_o$ as well as the fugacity equation for water since it can now be present in two phases. The ordering is very similar to the OG case; we keep the oil saturation as the third primary variable and move $y_3$ to the set of secondary variables. We previously mentioned that we need $y_3 = y_w$ in the set of secondary equations, and we now can be assured that the water conservation equation will never be an issue since water is present by construction in OWG cells. We order the secondary unknowns similarly to the OG case, making sure the fugacity constraints are properly aligned to their respective components.

\begin{figure}[t!]
    \centering
    \includegraphics[width=\textwidth]{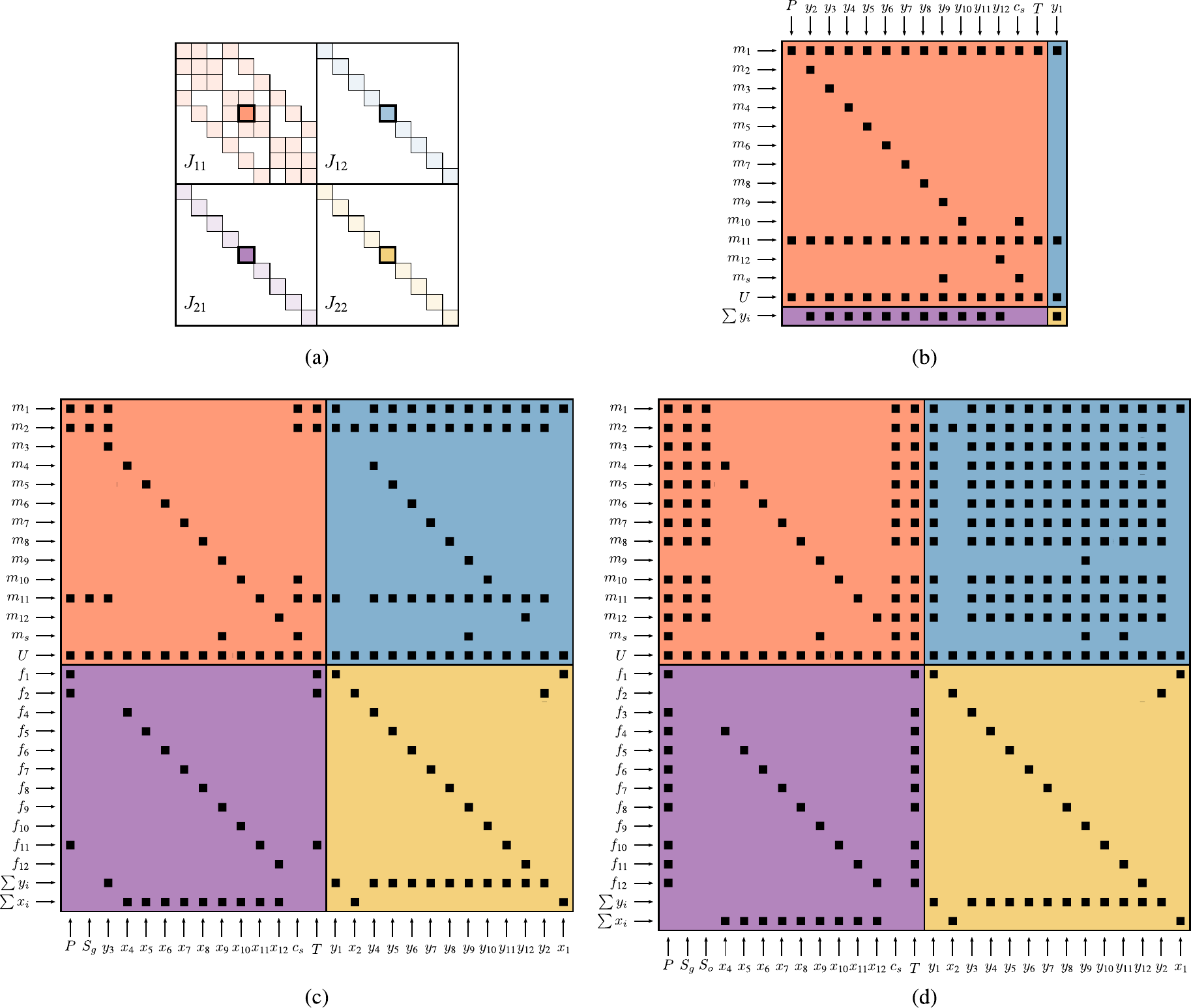}
    \caption{Jacobian system: (a) global partitioning into primary and secondary unknowns; (b) local Jacobian blocks for a gas (G) cell using the proposed ordering of unknowns; (c) same as (b) for an oil-gas (OG) cell; and (d) same as (b) for an oil-water-gas (OWG) cell. Each black square in (b), (c) and (d) is a non-zero value, denoting that the derivative of that equation with respect to the corresponding variable is non-zero. $m_i$ denotes the mass conservation equation of component $i$. The primary system is always $14\times14$ and the secondary systems are $1\times1$ (G), $13\times13$ (OG) and $14\times14$ (OWG).}
    \label{fig::ordering_final}
\end{figure}

\subsection{Summary}

In summary, we designed an ordering scheme for unknowns in the natural formulation to ensure that the Jacobian matrix is well-suited to the numerical methods we use to solve the linear system. Table \ref{tab::ordering} summarizes the set of unknowns in all three cases (G, OG, OWG), and Figure \ref{fig::ordering_final} shows the resulting sparsity pattern of diagonal Jacobian blocks.

\begin{table}[b!]
    \centering
    \caption{Primary and Secondary set of unknowns for all possible cases.}  
    {\footnotesize
    \begin{tabular}{lll}
    \toprule
    Case     & Flow Primary Unknowns & Flow Secondary Unknowns \\
    \midrule
    Gas & $[P,\;\!y_2,\;\!y_3,\;\!y_4,\dots,y_{n_c}]$ & $[y_1]$ \\
    Oil-Gas & $[P,S_g,\;\!y_3,x_4,\dots,x_{n_c}]$ & $[y_1,x_2,y_4,\dots,y_{n_c},y_2,x_1]$ \\
    Oil-Water-Gas & $[P,S_g,S_o,x_4,\dots,x_{n_c}]$ & $[y_1,x_2,y_3,y_4,\dots,y_{n_c},y_2,x_3]$ \\   
    \bottomrule
    \end{tabular}
    }
    \label{tab::ordering}
\end{table}
We note that although this ordering scheme has been designed for combustion cases, it will perform well for thermal-compositional cases (no reactions) as well. All of the previously mentioned considerations still apply since both phases saturations and components mole fractions can be driven to zero by displacement or vaporization.  With this particular partitioning, the reduced system $(\ref{eqn:reduced})$ may be readily assembled and solved with an appropriate preconditioning strategy.

\begin{remark}
Reservoir problems are typically driven by well boundary conditions, which have an important impact on the linear system and resulting solver strategy.  In this work, we eliminate all of the well unknowns using an exact Schur complement decomposition as suggested in \citet{Zhou13}. All of our cases use two wells with single perforations, allowing us to exactly invert the 2x2 block while preserving the sparsity pattern of the reservoir unknowns.  More sophisticated well treatments, while important, lie outside the scope of the current work.
\end{remark}

\section{Multi-Stage Preconditioning}
\label{sec::msp}

We now explore several preconditioning schemes for the primary system $A$.  These schemes differ in the specifics, but all of them may be cast as a split preconditioning,
\begin{equation}
    \left(N^{-1}AM^{-1}\right)\left(Mx\right) = N^{-1}b \,,
\end{equation}
with left-preconditioning operator $N^{-1}$ and a right-preconditioning operator $M^{-1}$.
The system matrix can be block factorized as $A=LDU$, where $D$ is a block diagonal matrix, and $L$ and $U$ are unit lower and upper block triangular matrices.  As a guiding principle, all of the preconditioners below attempt to approximate $N^{-1} \approx L^{-1}$ and $M^{-1} \approx (D U)^{-1}$. The details of how these approximations are made, however, will differ from one scheme to the next.

For all schemes, the left ``preconditioner'' is just a scaling operation that may be explicitly applied \emph{before} entering the Krylov solver.  We denote the scaled matrix and right-hand side as
\begin{equation} \label{eq::splitPrec1}
\bar A = N^{-1} A \,, \qquad \bar b = N^{-1} b \,.
\end{equation}
This cheap scaling pushes the system matrix closer to upper-block-triangular form.  After scaling, a multi-stage right-preconditioner is then applied within the Krylov iterations to solve the preconditioned system
\begin{equation} \label{eq::splitPrec2}
    \left(\bar A M^{-1}\right)\left(Mx\right) = \bar b \,.
\end{equation}
Multi-stage preconditioners are frequently used to efficiently tackle coupling in multi-physics problems.  The global form of a multi-stage preconditioner $M^{-1}$ can be formally written as,
\begin{equation} \label{eq::multistageprec}
M^{-1} = M_1^{-1} + \sum\limits_{i\,=\,2}^{n_\textrm{st}} M_i^{-1} \prod\limits_{j\,=\,1}^{i-1}\left(I - \bar AM_j^{-1}\right) \,,
\end{equation}
with $n_\textrm{st}$ the number of stages, $M_j^{-1}$ the $j^{th}$ preconditioner for $\bar A$, and $I$ the identity matrix.  In this work, we consider both two- and three-stage variants.

The system matrix contains three important groups of variables. It is therefore convenient to partition $A$ into a $3\times3$ block-system as
\begin{equation} 
A = \begin{bmatrix} A\sub{ss} & A\sub{sP} & A\sub{sT} \\ A\sub{Ps} & A\sub{PP} & A\sub{PT} \\ A\sub{Ts} & A\sub{Tp} & A\sub{TT} \end{bmatrix} \,,
\end{equation}
where the subscript $s$ denotes saturations and mole fractions, $P$ denotes pressure, and $T$ denotes temperature.  We then consider three preconditioner variants: 
\begin{enumerate}
\item A two-stage Constrained Pressure Residual (CPR) method
\item A two-stage Constrained Pressure Temperature Residual (CPTR) method
\item A three-stage Constrained Pressure Temperature Residual (CPTR3) method
\end{enumerate}
The three-stage preconditioner works directly on the $3\times3$ partitioning above.  In the two-stage variants, some variables are re-grouped to produce a $2\times2$ partitioning instead.
%
%

To unify the presentation, we adopt the following notational conventions: $A_{ij}$ indicates an individual block of the original system; $B_{ij}$ identifies a block of a first-level Schur complement approximation; $C_{ij}$ identifies a block of a second-level Schur complement approximation; $Z_{ij}$ is a block that may be well approximated with a zero block; and $D_{ij}$ denotes a (block-)diagonal matrix. 

\subsection{Constrained Pressure Residual (CPR)}

CPR leverages the different nature of the the mass conservation equations with respect to pressure (elliptic) and saturations (hyperbolic). Initially introduced in \citet{Wallis83,Wallis85} and extended in \citet{Lacroix03} and \citet{Cao05}, it was designed for isothermal, Black-Oil cases.  It is also used in compositional \citep{Voskov12b,Zhou13,Cusini15} and even thermal simulations \citep{Li15,Li17}, where it tends to show convergence issues. For thermal cases, it treats temperature as in the second stage.

To be precise, let the subscript $h$ to denote the union of saturation and temperature variables, with the re-partitioned system
\begin{equation}
A 
= \left[ \begin{array}{cc:c} A\sub{ss} & A\sub{sT} & A\sub{sP} \\ A\sub{Ts} & A\sub{TT} & A\sub{Tp} \\ \hdashline A\sub{Ps} & A\sub{PT} & A\sub{PP}\end{array} \right]
= \begin{bmatrix} A\sub{hh} & A\sub{hP} \\ A\sub{Ph} & A\sub{PP} \end{bmatrix} \,.
\end{equation}
Note that in the matrix $A_{hh}$ saturation and temperature variables are actually stored in an interleaved order, with small, dense blocks appearing for each grid cell.  The first step (left-scaling) used for CPR is:
\begin{equation} \label{eq::CPR}
\bar{A} = N^{-1}A = \begin{bmatrix} I & 0 \\ -D\sub{Ph}D\sub{hh}^{-1} & I \end{bmatrix} \begin{bmatrix} A\sub{hh} & A\sub{hP} \\ A\sub{Ph} & A\sub{PP} \end{bmatrix} = \begin{bmatrix} A\sub{hh} & A\sub{hP} \\ Z\sub{Ph} & B\sub{PP} \end{bmatrix} \,.
\end{equation}
Here, $D\sub{Ph}$ and $D\sub{hh}$ are block-diagonal approximations of $A\sub{Ph}$ and $A\sub{hh}$, respectively.  Also,
\begin{align}
Z\sub{Ph} & = A\sub{Ph} - D\sub{Ph}D\sub{hh}^{-1}A\sub{hh} \approx 0 \,, \\
B\sub{PP} & = A\sub{PP} - D\sub{Ph}D\sub{hh}^{-1}A\sub{hP}  \,.
\end{align}
If the block diagonal approximations were replaced with their exact equivalents, the effect of left-scaling would be to form an upper-block-triangular factor---i.e. $Z\sub{Ph}=0$.  The first stage of the right-preconditioner is then given by
\begin{equation} \label{eq::CPRM1}
M_1^{-1} = \begin{bmatrix} 0 & 0 \\ 0 & M\sub{PP}^{-1} \end{bmatrix}  \,,
\end{equation}
with $M\sub{PP}^{-1} = \texttt{amg}(B\sub{PP})$---i.e. the approximation of $B\sub{PP}^{-1}$ via an AMG preconditioner. We emphasize that $M^{-1}_1$ is not an invertible matrix, but instead a linear operator that approximates the action of the inverse of the pressure block. To correct errors associated with the as-yet-untouched saturation variables, the second stage is an ILU(0) sweep based on the full matrix.  It is given as,
\begin{equation}
M_2^{-1} = \texttt{ilu}\left(\bar{A}\right) \,.
\end{equation}

\begin{remark}
In the second-stage, it is often convenient to permute the matrix to interleave the pressure and hyperbolic variables, so that all unknowns are ordered cell-wise.  This interleaving is particularly useful if block-ILU (BILU) variants are chosen.
\end{remark}

\begin{remark}
In this work, we extract the diagonal of block matrices, such as $D\sub{Ph} = \texttt{diag}(A\sub{Ph})$, to use as approximations in the Schur-complements. The objective of the diagonal approximation is to preserve the sparsity pattern of the initial matrix in Schur-complement operations. These operators are similar to Quasi-IMPES operators \citep{Aziz79,Coa00}, but they include the energy equation in the decoupling for CPR.
\end{remark}

In thermal simulations, the energy equation can exhibit strongly elliptic behavior with respect to the temperature variable in the presence of high thermal diffusivity. In that case, treating it as if it were hyperbolic with respect to temperature in the second stage of CPR is unlikely to be a good approximation. In the next two preconditioner variants, we explore alternative treatments of the energy equation and temperature unknowns. 

\subsection{Two-stage Constrained Pressure-Temperature Residual (CPTR)}


We first consider a two-stage approach, where we now group temperature $T$ and pressure $P$ as elliptic variables, denoted with subscript $e$. This leads to the partitioning,
\begin{equation} \label{eq::CPTRinit}
A 
= \left[ \begin{array}{c:cc} A\sub{ss} & A\sub{sP} & A\sub{sT} \\ \hdashline A\sub{Ps} & A\sub{PP} & A\sub{PT} \\ A\sub{Ts} & A\sub{Tp} & A\sub{TT}\end{array} \right]
= \begin{bmatrix} A\sub{ss} & A\sub{se} \\ A\sub{es} & A\sub{ee} \end{bmatrix} \,.
\end{equation}
Similar to the CPR case, we use a lower-block triangular scaling to perform an approximate Schur reduction,
\begin{equation} \label{eq::CPTR1}
\bar{A} = N^{-1}A = \begin{bmatrix} I & 0 \\ -D\sub{es}D\sub{ss}^{-1} & I \end{bmatrix} \begin{bmatrix} A\sub{ss} & A\sub{se} \\ A\sub{es} & A\sub{ee} \end{bmatrix} = \begin{bmatrix} A\sub{ss} & A\sub{se} \\ Z\sub{es} & B\sub{ee} \end{bmatrix} \,,
\end{equation}
with
\begin{align}
Z\sub{es} & = A\sub{es} - D\sub{es}D\sub{ss}^{-1}A\sub{ss} \approx 0 \,, \\
B\sub{ee} & = A\sub{ee} - D\sub{es}D\sub{ss}^{-1}A\sub{se} \,.
\end{align}
%
%
To precondition this system, we need to approximate the inverse of the subsystem $B\sub{ee}$, which is a 2x2 block,
\begin{equation}
B\sub{ee}= \begin{bmatrix} B\sub{PP} & B\sub{PT} \\ B\sub{TP} & B\sub{TT} \end{bmatrix} \,.
\end{equation}
One option is to use a dedicated block preconditioner, similar to the algorithm in \citet{Roy19a}, but due to the complexity of the compositional formulation here we could not pursue a similar strategy. To tackle this block monolithically, we would need a well-suited algebraic preconditioner designed for systems of PDEs. Example options include System-AMG (SAMG) \citep{Gries15} or BoomerAMG \citep{Baker11} from the Hypre library \citep{Falgout02}. The latter has shown promising results for coupled PDEs arising from multiphase flow \citep{Bui17,Bui18}, including two-phase dead-oil thermal water injection \citep{Roy19b} and multiphase poromechanics \citep{Bui20}.

The first-stage preconditioner is then simply
\begin{equation} \label{eq::CPTRM1}
M_1^{-1} = \begin{bmatrix} 0 & 0 \\ 0 & M\sub{ee}^{-1} \end{bmatrix} \,.
\end{equation}
The second-stage preconditioner is identical to CPR,
\begin{equation}
M_2^{-1} = \texttt{ilu}\left(\bar{A}\right) \,.
\end{equation}
This two-stage variant will be denoted CPTR in the remainder of this work. 


\subsection{Three-stage Constrained Pressure-Temperature Residual (CPTR3)}
Finally, we consider a three-stage variant.  In the previous approach, the first left scaling step leads to
\begin{equation} \label{eq::CPTR31}
\hat{A} = N_1^{-1}A = \begin{bmatrix} I & 0 & 0 \\ -D\sub{Ps}D\sub{ss}^{-1} & I & 0 \\ -D\sub{Ts}D\sub{ss}^{-1} & 0 & I \end{bmatrix} \begin{bmatrix} A\sub{ss} & A\sub{sP} & A\sub{sT} \\ A\sub{Ps} & A\sub{PP} & A\sub{PT} \\ A\sub{Ts} & A\sub{Tp} & A\sub{TT} \end{bmatrix} = \begin{bmatrix} A\sub{ss} & A\sub{sP} & A\sub{sT} \\ Z\sub{Ps} & B\sub{PP} & B\sub{PT} \\ Z\sub{Ts} & B\sub{Tp} & B\sub{TT} \end{bmatrix} \,,
\end{equation}
with
\begin{align}
Z\sub{Ps} & = A\sub{Ps} - D\sub{Ps}D\sub{ss}^{-1}A\sub{ss} \approx 0 \\
Z\sub{Ts} & = A\sub{Ts} - D\sub{Ts}D\sub{ss}^{-1}A\sub{ss} \approx 0 \\
B\sub{PP} & = A\sub{PP} - D\sub{Ps}D\sub{ss}^{-1}A\sub{sP} \,, \\
B\sub{PT} & = A\sub{PT} - D\sub{Ps}D\sub{ss}^{-1}A\sub{sT} \,, \\
B\sub{Tp} & = A\sub{Tp} - D\sub{Ts}D\sub{ss}^{-1}A\sub{sP} \,, \\
B\sub{TT} & = A\sub{TT} - D\sub{Ts}D\sub{ss}^{-1}A\sub{sT} \,.
\end{align}
To get closer to an upper block-triangular system, we perform a second scaling:
\begin{equation} \label{eq::CPTR32}
\bar{A} = N_2^{-1}\hat{A} = \begin{bmatrix} I & 0 & 0 \\ 0 & I & 0 \\ 0 & -D\sub{Tp}D\sub{PP}^{-1} & I \end{bmatrix} \begin{bmatrix} A\sub{ss} & A\sub{sP} & A\sub{sT} \\ Z\sub{Ps} & B\sub{PP} & B\sub{PT} \\ Z\sub{Ts} & B\sub{Tp} & B\sub{TT} \end{bmatrix} = \begin{bmatrix} A\sub{ss} & A\sub{sP} & A\sub{sT} \\ Z\sub{Ps} & B\sub{PP} & B\sub{PT} \\ Z\sub{Ts} & Z\sub{Tp} & C\sub{TT} \end{bmatrix} \,,
\end{equation}
with
\begin{align}
Z\sub{Ts} &= Z\sub{Ts} - D\sub{Tp}D\sub{PP}^{-1}Z\sub{Ps} \approx 0 \,, \\
Z\sub{Tp} &= Z\sub{Tp} - D\sub{Tp}D\sub{PP}^{-1}B\sub{PP} \approx 0 \,, \\
C\sub{TT} &= B\sub{TT} - D\sub{Tp}D\sub{PP}^{-1}B\sub{PT} \,.
\end{align}
The first stage preconditioner is then 
\begin{equation}\label{eq::CPTR3M1}
M_1^{-1} = \begin{bmatrix} 0 & 0 & 0 \\ 0 & 0 & 0 \\ 0 & 0 & M\sub{TT}^{-1} \end{bmatrix} \,,
\end{equation}
with $M\sub{TT}^{-1} = \texttt{amg}\left(C_{TT}\right)$. The second stage is
\begin{equation} \label{eq::CPTRM2}
M_2^{-1} = \begin{bmatrix} 0 & 0 & 0 \\ 0 & M\sub{PP}^{-1} & 0 \\ 0 & 0 & 0 \end{bmatrix} \,,
\end{equation}
with $M\sub{PP}^{-1} = \texttt{amg}\left(B\sub{PP}\right)$. The last stage is the local ILU(0) sweep of the scaled matrix,
\begin{equation}\label{eq::CPTR3M3}
M^{-1}_3 = \texttt{ilu}\left(\bar{A}\right) \,.
\end{equation}
This three-stage variant will be denoted CPTR3 in the remainder of this work. Compared to CPR, we have one more subsystem solve for the temperature matrix. It is important to note that this algorithm can be readily added on top of any existing CPR implementation. It only requires block matrices and scalar AMG preconditioners, both prerequisites for CPR.

\subsection{Summary Algorithms}

The procedure to apply the two-stage preconditioners to a vector $(y = M^{-1}r)$ are given in Algorithm \ref{alg::two}, and the procedure for the three-stage preconditioner in Algorithm \ref{alg::three}.
\begin{algorithm}[H]
{\footnotesize
\caption{Apply two-stage preconditioner}
\label{alg::two}
\begin{algorithmic}[1]
\Function{ $y$ = \tt applyTwoStagePrec}{$r$}
    \State $x = M^{-1}_1r$
    \State $z = r - \bar Ax$
    \State $y = x + M_2^{-1}z$
	\State \Return{$y$}
 \EndFunction
\end{algorithmic}
}
\end{algorithm}
\vspace{-1.25em}
\begin{algorithm}[H]
{\footnotesize
\caption{Apply three-stage preconditioner}
\label{alg::three}
\begin{algorithmic}[1]
\Function{ $y$ = \tt applyThreeStagePrec}{$r$}
    \State $x = M^{-1}_1r$
    \State $z = r - \bar Ax$
    \State $v = M^{-1}_2z$  
    \State $u = z - \bar Av$
    \State $y = x + v + M_3^{-1}u$
	\State \Return{$y$}
 \EndFunction
\end{algorithmic}
}
\end{algorithm}

\section{Numerical Results}
\label{sec::results}

\subsection{Thermal-Compositional Homogeneous Cases}

We start by studying thermal-compositional cases, with no reactions and homogeneous properties (porosity and permeability). For all cases, we run the simulator (AD-GPRS) on a laboratory scale, 2D case using parameters for an extra-heavy oil from Venezuela \citep{Lapene10b}. Our compositional description has seven hydrocarbon components, nitrogen, oxygen, carbon dioxide, water and one coke solid species. This gives $n_c = 12$ and $n_s=1$, for a total of 14 unknowns per cell. Table \ref{tab::comp-therm} shows the main parameters of our simulations. We inject hot air into a mixture of oil, water and nitrogen for 100 minutes. The water phase gets vaporized by the hot temperature front, and light hydrocarbon components will be stripped from the oil phase. Although no chemical reactions are taking place, this case shows a strong coupling between the mass and energy transport through viscosity, density and relative permeability calculations.

\begin{table}[htb]
    \centering
    \caption{Parameters for compositional--thermal simulations.}
    {\footnotesize
    \begin{tabular}{llrl}
    \toprule
    Property & Symbol & Value & Unit \\
    \midrule
    Domain Size & $L$ & 0.35 & m \\
    Porosity & $\phi$ & 0.36 & -- \\
    Permeability & $k$ & 10 & D \\
    Injection Rate & $q$ & 4.32 & m$^3$/day \\
    Injection Temperature & $T_\textrm{inj}$ & 873.15 & K \\
    Initial Temperature & $T_\textrm{init}$ & 323.15 & K \\
    Initial Pressure & $P_\textrm{init}$ & 7.8 & bar \\
    Initial Oil Saturation & $S_o$ & 0.4791 & -- \\
    Initial Water Saturation & $S_w$ & 0.2048 &  -- \\
    \bottomrule
    \end{tabular}
    }
    \label{tab::comp-therm}
\end{table}

As representative linear test problems, we have output the Jacobian matrix and residual vector from the simulator at multiple time-steps at the first Newton iteration. All of the convergence results are then obtained with a right-preconditioned GMRES \citep{Saad86} algorithm applied to the left-scaled system given in Equation \eqref{eq::splitPrec1}, using a relative tolerance of $10^{-8}$ and no restart. 
Using homogeneous properties allows us to conduct a mesh refinement study and quantify the impact of the grid size on the linear solver convergence.

The absence of an enthalpy source (from reactions) in the domain allows us to quantify the relative magnitude of thermal advection versus diffusion in a straightforward way. The effects can be compared through the dimensionless thermal P\'eclet number \citep{Incropera07},
\begin{equation} \label{eq::pe}
    \textrm{Pe} = \dfrac{uL\rho c_p}{\kappa} = \dfrac{q\rho c_p}{\kappa L} \,,
\end{equation}
with $L$ a characteristic length, $u$ a characteristic velocity, $q$ a characteristic flow rate, $\rho$ and $c_p$  the density and the specific heat capacity of the fluid carrier and $\kappa$ the global thermal conductivity. With no conduction ($\kappa=0$), Pe is infinite, and with no convection ($c_p = 0$) Pe is zero. The only heat source in the absence of reactions is  the injection well, where we know the composition, flow rate, and temperature of the fluid. We use this data to compute the density and heat capacity of the air for  the P\'eclet number computation, as well as the air injection flow rate. The characteristic length is the length of the domain, and the thermal conductivity an input parameter depending on the rock properties. Table \ref{table::Pe} summarizes the values we use to compute our P\'eclet numbers. Note that to vary Pe, we will use the thermal conductivity of the rock ($\kappa$) while keeping all fluid properties constant.

\begin{table}[htb]
    \centering
    \caption{Parameters for the P\'eclet number calculations, resulting in Pe = $\mathcal{O}(1)$. Properties for air are taken at injection conditions: 600$^\circ$C and $\sim$8 bars.}
    {\footnotesize
    \begin{tabular}{crl}
        \toprule
        Parameter & Value & Unit \\
        \midrule
        $q$ & 4.32 & m$^3$/day \\
        $L$ & 0.35 & m \\
        $c_p$ & 1.42 & kJ/kg/K \\
        $\kappa$ & 50 & kJ/m/day/K \\
        $\rho$ & 3.1 & kg/m$^3$ \\
        \bottomrule
    \end{tabular}
    }
    \label{table::Pe}
\end{table}

\begin{figure}[t]
    \centering
    \includegraphics[width=.7\textwidth]{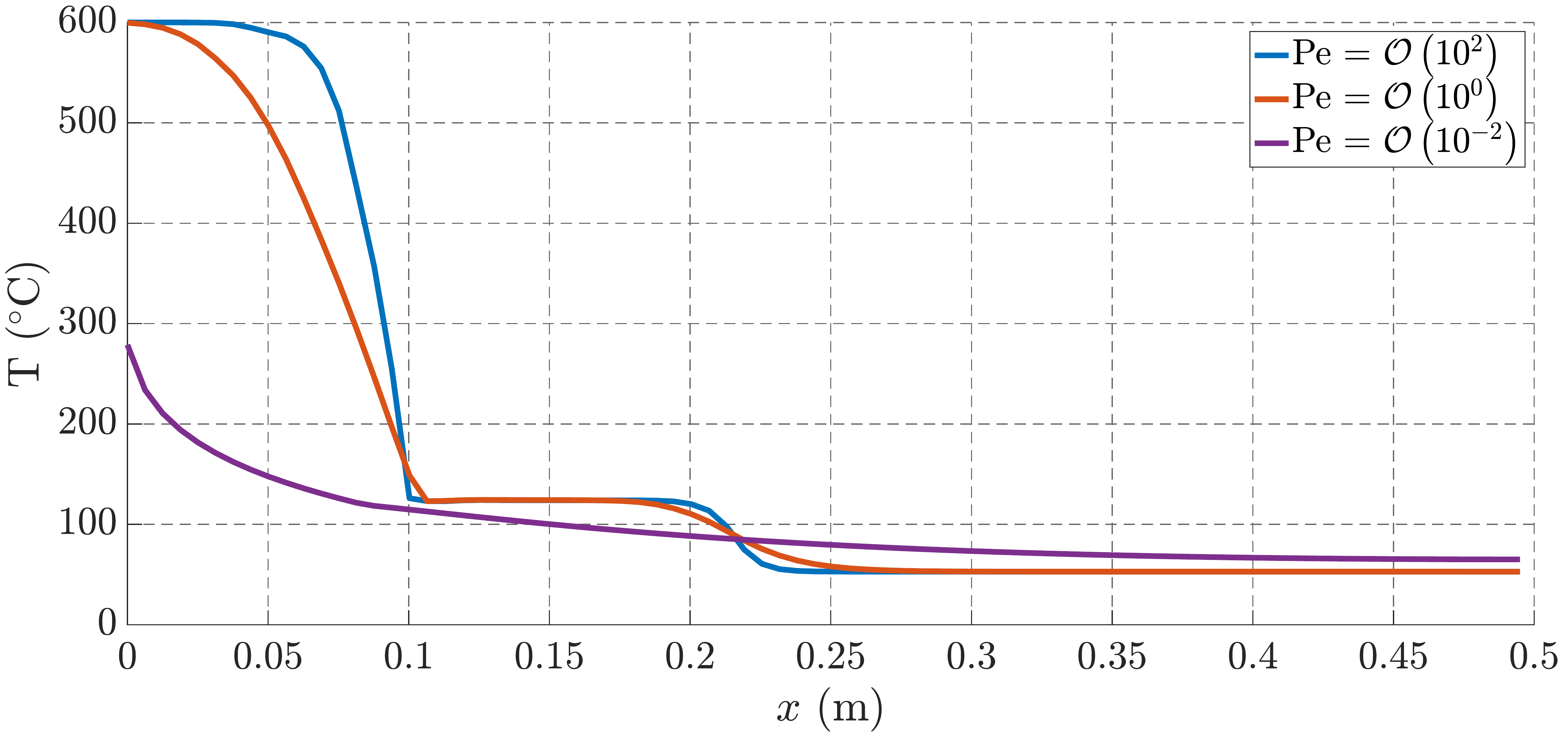}
    \caption{Temperature cross-section between wells for different Pe numbers at $t=50$ min.}
    \label{fig::Pe}
\end{figure}

Figure \ref{fig::Pe} plots the temperature profile for three different P\'eclet numbers along a cross-section through the domain. We can clearly observe different regimes, transitioning from advection dominated with sharp temperature fronts (blue curve) to a diffusion dominated case with a smooth profile and clear long-range interactions (purple curve). The orange curve is a unit P\'eclet number case, in the transition regime. The nature of the energy equation will change according to the P\'eclet number, from hyperbolic for pure advection cases to elliptic for pure diffusion cases, with  intermediate cases being parabolic. 

\begin{remark} We present numerical results here modeling laboratory scale specimens, because our reaction model is only valid for small grid block sizes (see section \ref{subsec::isc}) and all cases are based on reactive ones. However, since this study is based on the dimensionless version of the energy equation, the results will also hold at the reservoir scale, and the range of P\'eclet number we study is very representative of steam/gas injection field cases \citep{Prats82}.
\end{remark}

We start by studying the quality of the various Schur complement approximations themselves. To do so, we use a direct solver (rather than AMG) to apply the \emph{exact} operators,
\begin{align}
M\sub{PP}^{-1} &= B\sub{PP}^{-1}\ \ \ \ (\textrm{CPR and CPTR3}) \,,\\
M\sub{TT}^{-1} &= C\sub{TT}^{-1}\ \ \ \ (\textrm{CPTR3}) \,,\\
M\sub{ee}^{-1} &= B\sub{ee}^{-1}\ \ \ \ (\textrm{CPTR}) \,.
\end{align}

\begin{figure}[b!]
    \centering
    \includegraphics[width=\textwidth]{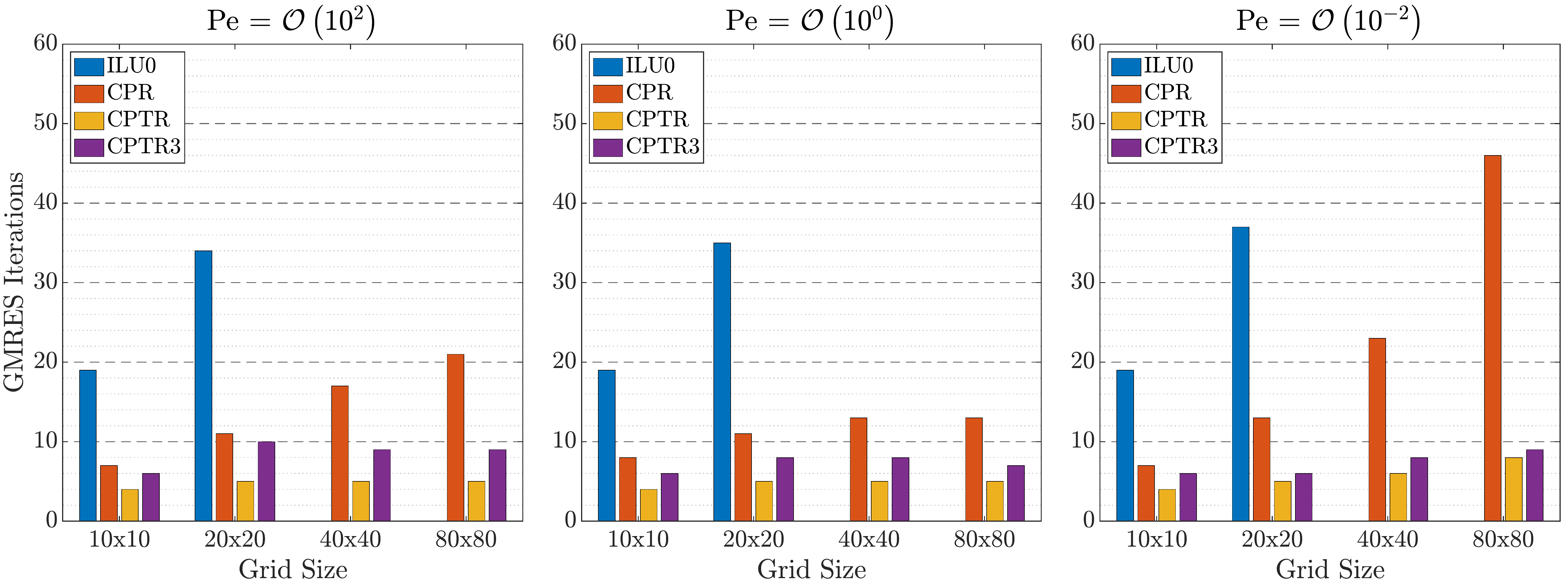}
    \caption{GMRES iterations results for the advection dominated (left), neutral (center) and diffusion dominated (right) cases using a direct solver for all subsystems.}
    \label{fig::results25}
\end{figure}

Using a direct solver also allows for a straightforward implementation of CPTR, since we do not have to worry about a multi-PDE AMG routine. Figure \ref{fig::results25} summarizes results for the $t=25$ min time step, and a complete performance profile is given in Table \ref{tab::direct25min}. Note that the results correspond to the first non-linear iteration of that time step. Unless otherwise specified, the time step is constant and set at 10 seconds.

\begin{table}[t!]
    \centering
    \caption{Thermal--Compositional results for the first non-linear iteration at $t = 25$ min. An asterisk ($^*$) denotes cases showing grid orientation numerical issues, unrelated to the linear solver.}
    {\footnotesize
    \begin{tabular}{lcrrcp{0.7cm}p{0.7cm}p{0.7cm}p{0.8cm}p{1.3cm}p{1.6cm}}
    \toprule
     & & & & & \multicolumn{6}{c}{GMRES Iterations} \\
    \cmidrule(){6-11}
    P\'eclet & Grid Size & Matrix Size & Non-zeros & $\ $ & ILU0 & CPR & CPTR & CPTR3 & CPR-AMG & CPTR3-AMG \\
    \midrule
    \multirow{6}{*}{$\mathcal{O}\left(10^{2}\right)$} & 10x10 & 1,400 & 41,000 & & 19 & 7 & 4 & 6 & 5 & 7 \\
     & 20x20 & 5,600 & 172,000 & & 34 & 11 & 5 & 10 & 11 & 10 \\
     & 40x40 & 22,400 & 694,000 & & 73 & 17 & 5 & 9 & 27 & 12 \\
     & 80x80 & 89,600 & 2,790,000 & & 142 & 21 & 5 & 9 & 42 & 15 \\
     & 125x125 & 218,750 & 6,904,000 & & 207$^*$ & 43$^*$ & 8$^*$ & 12$^*$ & 104$^*$ & 37$^*$ \\
     & 160x160 & 358,400 & 11,328,000 & & 277$^*$ & 46$^*$ & 9$^*$ & 22$^*$ & 190$^*$ & 42$^*$ \\
    \midrule
    \multirow{6}{*}{$\mathcal{O}\left(10^{0}\right)$} & 10x10 & 1,400 & 42,000 & & 19 & 8 & 4 & 6 & 8 & 6 \\
     & 20x20 & 5,600 & 173,000 & & 35 & 11 & 5 & 8 & 11 & 9 \\
     & 40x40 & 22,400 & 697,000 & & 70 & 13 & 5 & 8 & 16 & 10 \\
     & 80x80 & 89,600 & 2,802,000 & & 142 & 13 & 5 & 7 & 20 & 12 \\
     & 125x125 & 218,750 & 6,930,000 & & 219 & 14 & 5 & 7 & 22 & 14 \\
     & 160x160 & 358,400 & 11,226,000 & & 281 & 16 & 6 & 9 & 26 & 16 \\
    \midrule
    \multirow{6}{*}{$\mathcal{O}\left(10^{-2}\right)$} & 10x10 & 1,400 & 42,000 & & 19 & 7 & 4 & 6 & 8 & 6 \\
     & 20x20 & 5,600 & 173,000 & & 37 & 13 & 5 & 6 & 13 & 7 \\
     & 40x40 & 22,400 & 705,000 & & 71 & 23 & 6 & 8 & 24 & 10 \\
     & 80x80 & 89,600 & 2,846,000 & & 138 & 46 & 8 & 9 & 50 & 15 \\
     & 125x125 & 218,750 & 6,964,000 & & 210 & 80 & 9 & 14 & 89 & 21 \\
     & 160x160 & 358,400 & 11,423,000 & & 251$^*$ & 118$^*$ & 12$^*$ & 23$^*$ & 136$^*$ & 33$^*$ \\
    \bottomrule
    \end{tabular}
    }
    \label{tab::direct25min}
\end{table}

We compare the multi-stage preconditioners with a single-stage ILU(0) preconditioner. The latter approach does not take any physics into account nor does it scale well for this problem. Unless the grid size is very small, it will not converge in a reasonable number of iterations for any of our cases. We remark that the finer ($160\times160$) grid size exhibits strong grid-orientation effects when the P\'eclet number is far from unity, and therefore linear solver performance trends are likely corrupted for these cases. The reader is referred to \citet{Kozdon09} for more details about grid-orientation problems for transport in porous media.  While important, these challenges are outside the scope of the current work. 

With respect to convergence behavior for increasing grid resolution, we observe that both CPTR and CPTR3 are virtually unaffected. As expected, CPR struggles to converge when the P\'eclet number is low, since the conduction part of the temperature dominates and the ILU(0) second stage cannot efficiently reduce the low frequency modes in the energy equation. Those convergence issues worsen with the grid size. 

CPTR here solves the pressure/temperature subsystem (of size 2$n_b$) exactly using a direct solver. One then observes excellent iteration performance---though not timing performance due to the cost of the direct solves---regardless of the P\'eclet number and the grid size. In this work, the focus in on the algorithmic behavior of the methods, hence we chose to report iteration counts as the metric of interest to assess the merit of each approach. Performance studies based on CPU times will be the subject of future work. The ILU(0) second stage is able to efficiently reduce the remaining high frequency error modes in the saturation/mole fraction variables. For the $80\times80$ grid, CPTR outperforms CPR in terms of iteration counts by 76\%, 48\% and 83\% for high, unit and low P\'eclet numbers respectively. CPTR3 introduces one more level of approximation than CPTR, and as a result it shows a slightly higher number of iterations. However, it is still able to perform well and for an overall lower computational cost. For the $80\times80$ grid, CPTR3 outperforms CPR in terms of iteration counts by 57\%, 46\% and 80\% for high, unit and low P\'eclet numbers respectively. These results demonstrate that a specific treatment of the energy equation and temperature unknown coupling significantly increases the performance and robustness of the preconditioners. The performance of the proposed methods is always better than CPR across all grid sizes and thermal regimes. 



\begin{figure}[t]
    \centering
    \includegraphics[width=\textwidth]{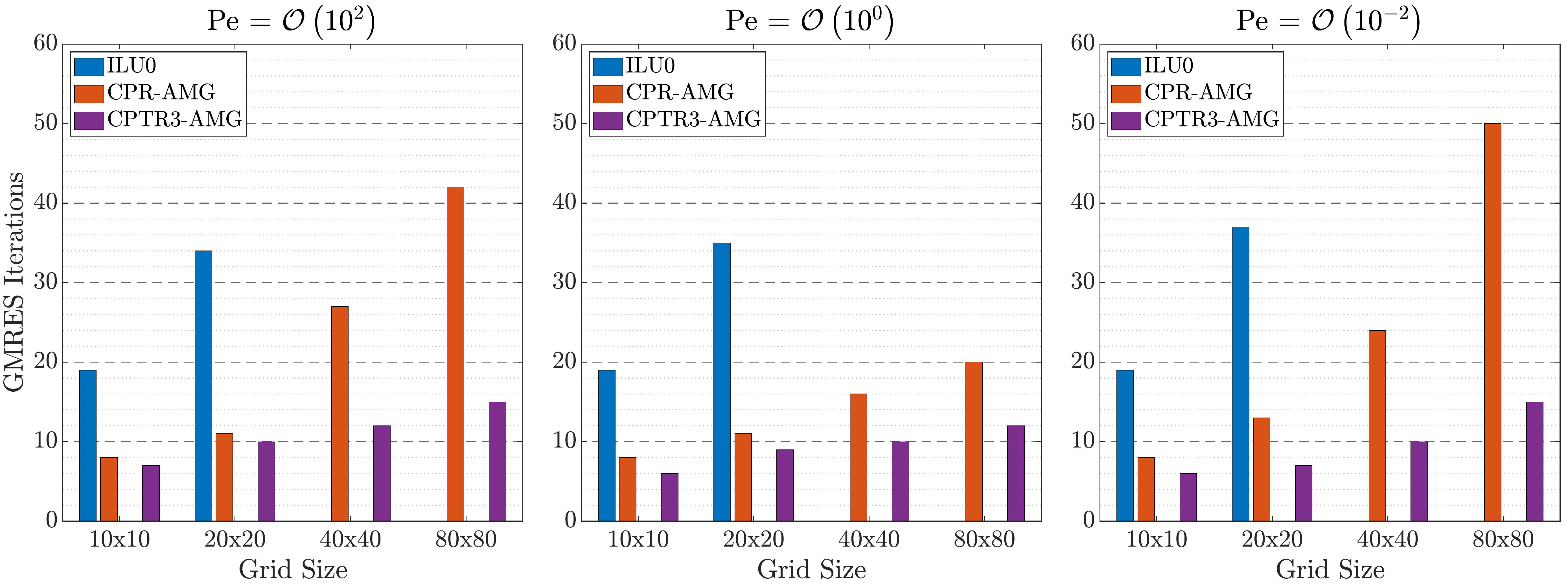}
    \caption{GMRES iterations results for the advection dominated (left), neutral (center) and diffusion dominated (right) cases using AMG preconditioners for all subsystems.}
    \label{fig::resultsAMG25}
\end{figure}

We now consider a more scalable implementation of the preconditioners, using AMG for the key subsystems.  We denote the AMG versions of the preconditioners by adding the suffix ``-AMG". Figure \ref{fig::resultsAMG25} shows the results for $t = 25$ min, and a complete performance summary is also given in Table \ref{tab::direct25min}. Given the good results using CPTR3 in the previous section, and the fact it only requires a scalar AMG implementation, we focus on this approach in the remainder of the work. CPTR remains an appealing strategy, and could perhaps outperform CPTR3, but it requires a robust system-AMG implementation. For our comparison, we use a classic AMG method \citep{Ruge87} as implemented in the \texttt{HSL\_MI20} package \citep{hsl} with default parameters: symmetric Gauss-Seidel smoother, single V-cycle, and a direct coarse solver. We set the maximum coarse problem size to 1000 degrees of freedom.

CPTR3-AMG is able to keep the number of GMRES iterations below 15 in all cases up to the $80\times80$ grid size. Both the pressure and temperature subsystems are well approximated and the scaling with respect to grid size is good. Using the $80\times80$ grid, CPTR3-AMG outperforms CPR-AMG by 64\%, 36\% and 74\% for high, unit and low P\'eclet numbers respectively. CPR-AMG shows an increased number of iterations compared to CPR for unit and low P\'eclet numbers of about 50\%. However, for the high P\'eclet number case, we see a much larger increase, leading to 150\% more iterations compared to the direct sub-solver version.

\begin{figure}[t!]
    \centering
    \includegraphics[width=.85\textwidth]{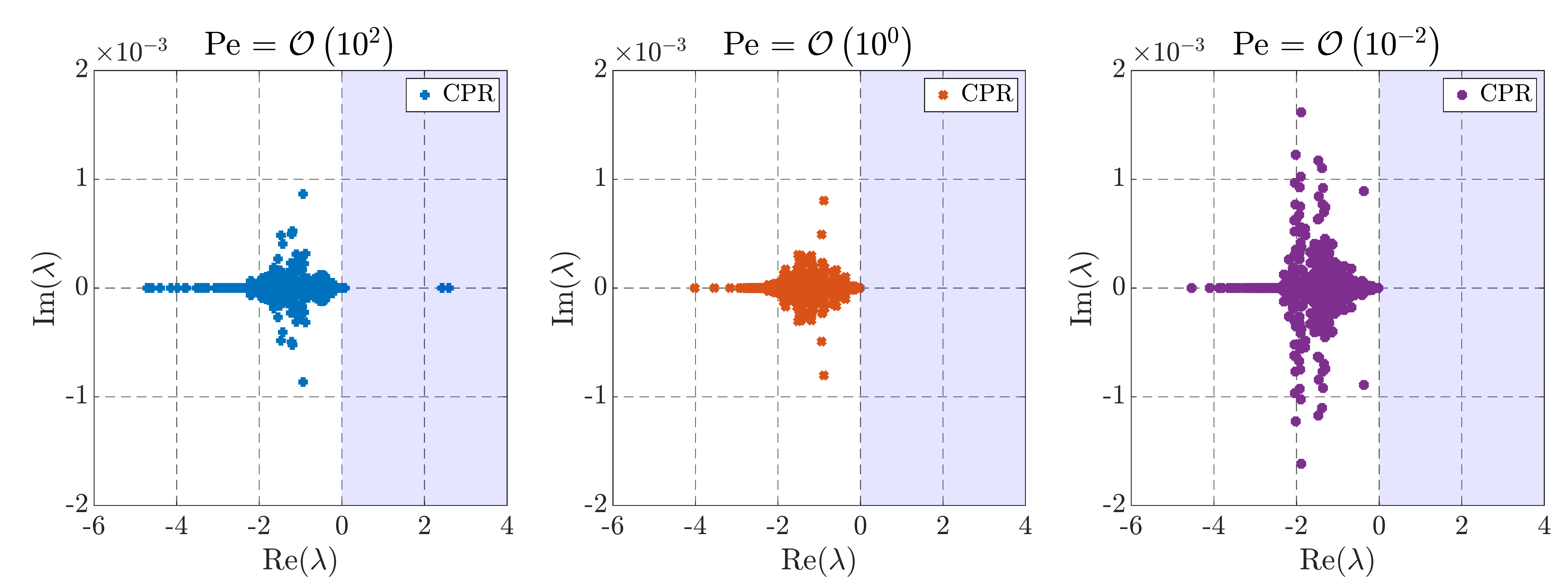}
    \includegraphics[width=.85\textwidth]{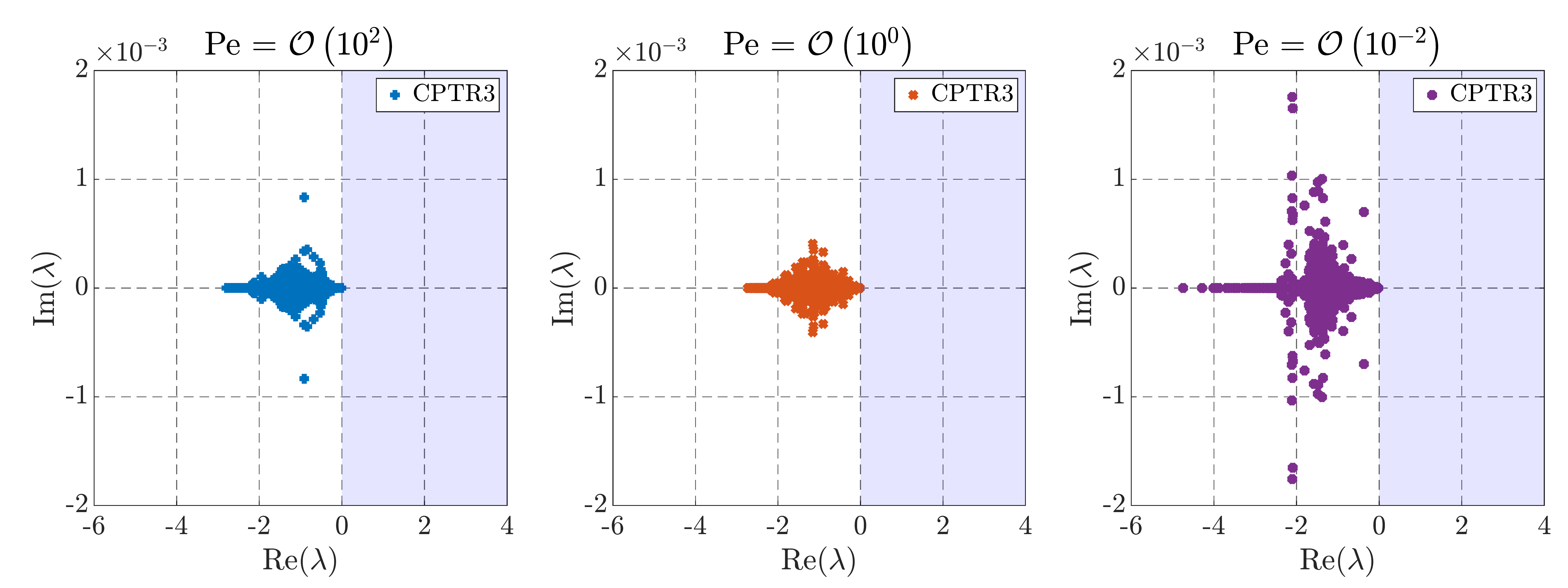}
    \caption{Eigenvalues of pressure system for the CPR method (top three) and the CPTR3 method (bottom three). In both cases, we show the results for the high (left), unit (middle) and low (right) P\'eclet number, computed for the 80x80 grid.}
    \label{fig::eigsApp}
\end{figure}

To investigate further, we compute the eigenvalues of the pressure matrix $B\sub{PP}$ for the $80\times80$ grid and the CPR preconditioner. Figure \ref{fig::eigsApp} (top) shows the spectra for the high (left), unit (middle) and low (right) P\'eclet numbers using the CPR preconditioner. For the high P\'eclet number we find both positive and negative eigenvalues. In the CPR case, the first pressure decoupling step includes the energy equation, since the temperature is part of the $h$ variables (Eq. \ref{eq::CPR}). It seems doing so corrupts the ellipticity of the $B\sub{PP}$ matrix. We observe better behavior with the $B_{PP}$ matrix we get with the CPTR3 preconditioner, for which the pressure decoupling is done purely with mole conservation equations. That matrix is negative definite as desired, as shown in Figure \ref{fig::eigsApp} (bottom).
 
\begin{figure}[b!]
    \centering
    \includegraphics[width=\textwidth]{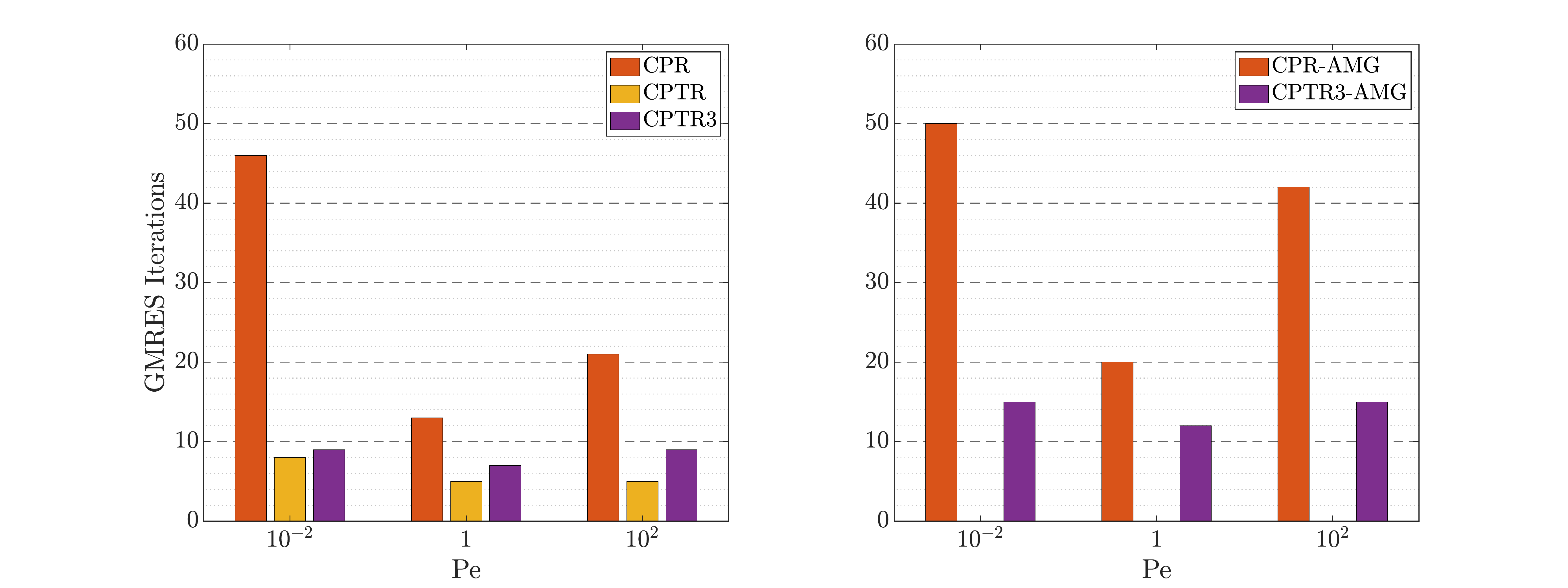}
    \caption{Number of GMRES iterations for the $80\times80$ grid for multiple P\'eclet numbers, using direct sub-solvers (left) and AMG sub-preconditioners (right).}
    \label{fig::MultPeAMG}
\end{figure}

A very important feature for a preconditioner is its ability to deal with different physical regimes that may arise over the course of the simulation. In thermal reservoir simulations, different parts of a reservoir can show different property values and local P\'eclet numbers can greatly vary. If the reservoir is composed of different rock types, the thermal conductivity can also vary and impact the local P\'eclet number. Figure \ref{fig::MultPeAMG} illustrates that the sensitivity to  thermal regime is greatly reduced for the proposed preconditioners, both using direct solvers (top), and using AMG preconditioners (bottom). To provide a quantitative metric, we compute the coefficient of variation (CV) of the GMRES iterations across P\'eclet numbers, defined as the standard deviation divided by the mean. For CPR-AMG, we get a CV of 41.6\%, but that number reduces to 12.8\% for CPTR3-AMG.

\subsection{Thermal--Compositional Heterogeneous Cases}

So far we have considered lab scale cases with homogeneous properties (leading to a fairly uniform velocity profile), but in a typical reservoir the fluid velocities can exhibit differences of several orders of magnitude.

\begin{figure}[t]
    \centering
    \includegraphics[width=0.7\textwidth]{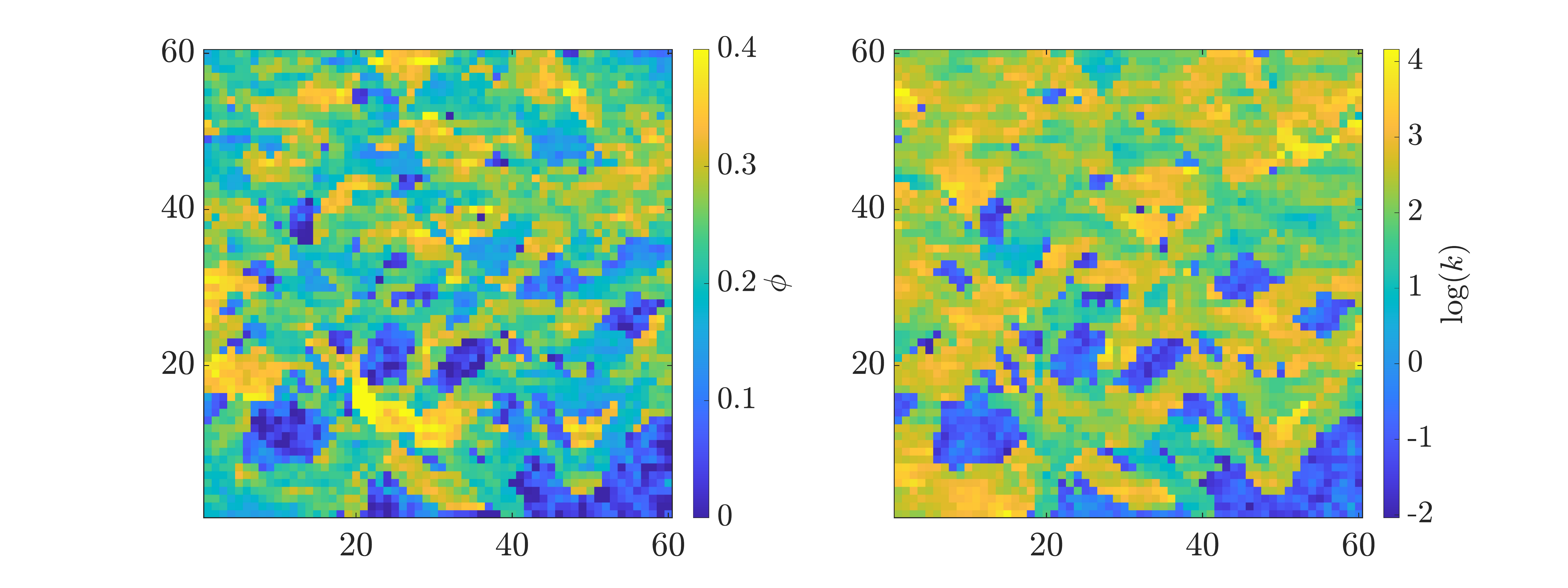}
    \includegraphics[width=0.7\textwidth]{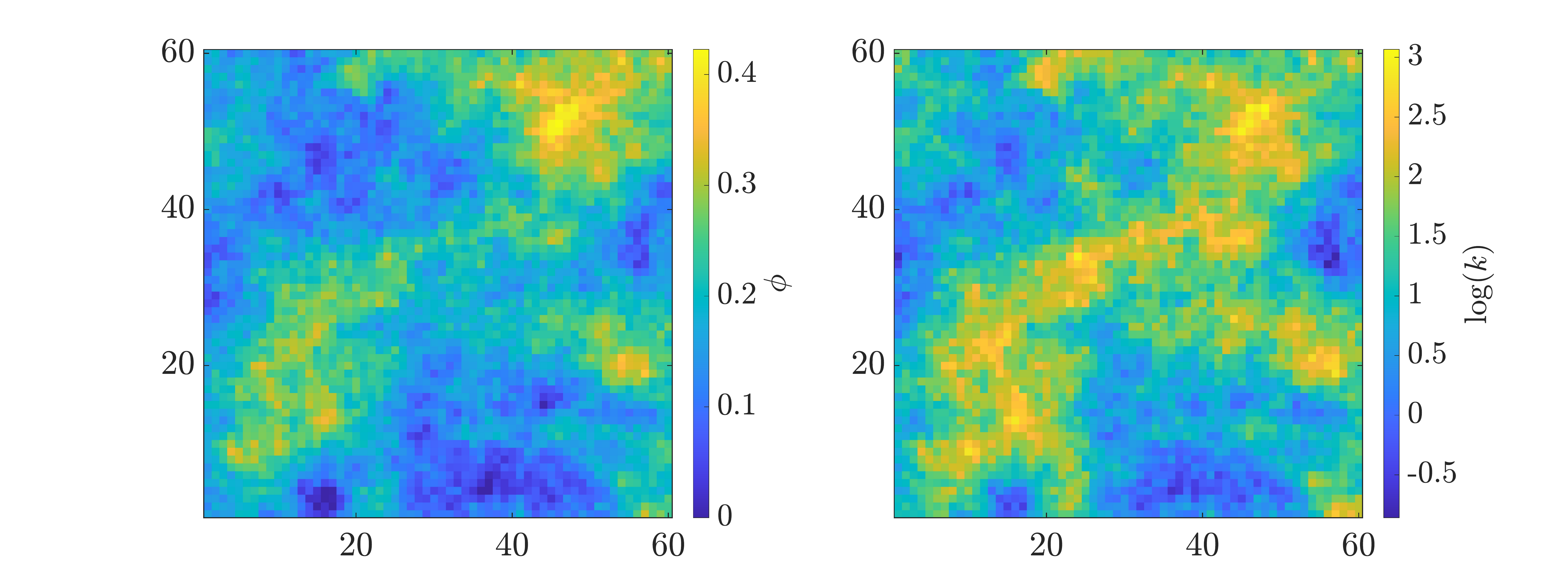}
    \caption{Porosity (left) and log-permeability (right) from the top layer (top) and bottom layer (bottom) of SPE10 ($60\times60$). They model a prograding environment with smooth property variations and a fluvial environment with sharp property variations, respectively.}
    \label{fig::SPE10}
\end{figure} 

We now test our method on heterogeneous cases to confirm the trends and performance we observed in the previous section. We generate two different problems based on the SPE10 Model 2 dataset \citep{Christie01}. These cases are identical to the homogeneous cases for all parameters and properties except for the porosity and permeability fields. We crop the initial SPE10 layers ($220\times60$ cells) to 60 cells in the x-direction to create two square patterns. We retain a quarter-five spot pattern, injecting at the bottom left corner and producing at the top right corner. Both squares are 0.36 meters wide and show a similar length scale and pore volume compared the homogeneous cases. We plot the porosity and log-permeability fields on Figure \ref{fig::SPE10}. Note that we do not vary the grid refinement here, in order to avoid any issues with property downscaling for a consistent comparison.

\begin{figure}[b!]
    \centering
    \includegraphics[width=\textwidth]{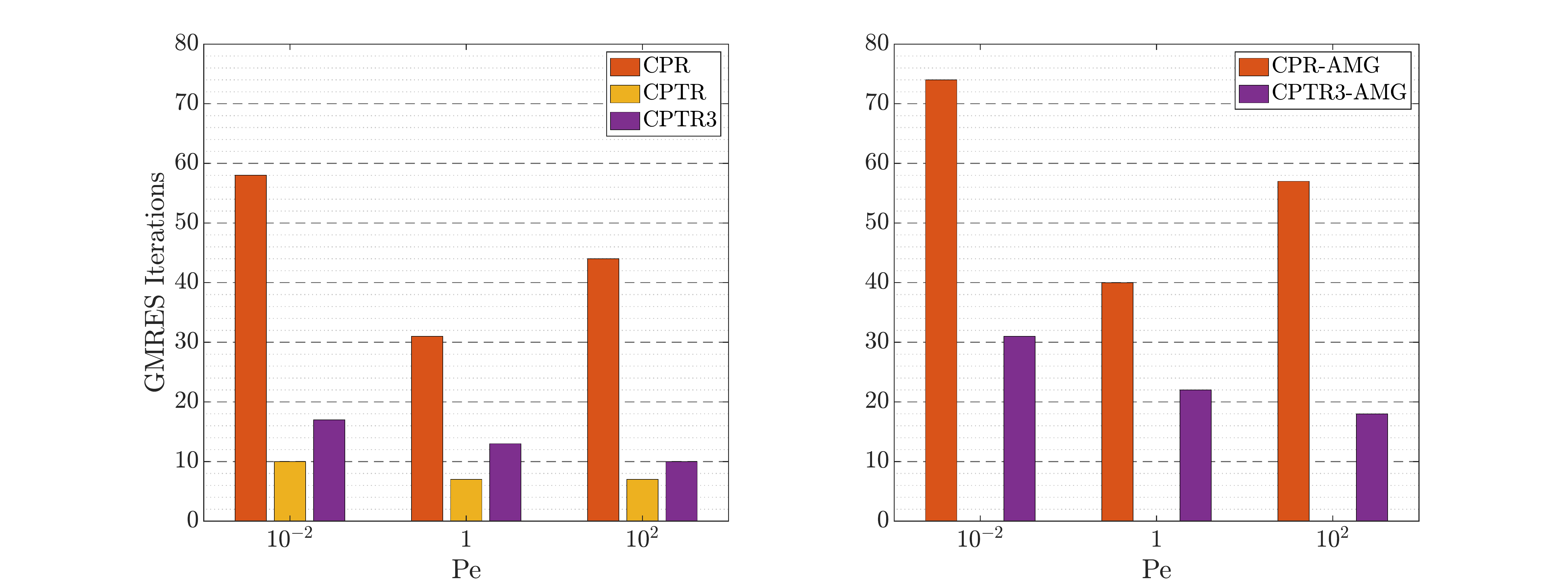}
    \caption{SPE10 top layer, results at $t=15$ min, using direct sub-solvers (left) and AMG sub-preconditioners (right).}
    \label{fig::MultPeTop}
\end{figure}

We start with the top layer, modeling a prograding near-shore environment. The property variations are smooth but the heterogeneity is quite strong, with the permeability ranging across four orders of magnitude. The problem has 50,400 primary unknowns with $\sim$1.6 millions non-zeros in the system matrix. We conduct the same numerical study as for the homogeneous case. Figure \ref{fig::MultPeTop} shows the performance of both new preconditioners using direct solvers for the subsystems (left) and using AMG preconditioners (right). These results confirm the trend we observed in the homogeneous case, and we actually see even more improvement. The addition of heterogeneity strengthens the coupling between transport and thermal effects, since now the viscosity reduction is spatially dependent. For all P\'eclet numbers, the number of iterations remains under 31 with an average of 24 using CPTR3-AMG, whereas CPR-AMG requires at least 40 iterations and an average of 57.

\begin{figure}[t]
    \centering
    \includegraphics[width=\textwidth]{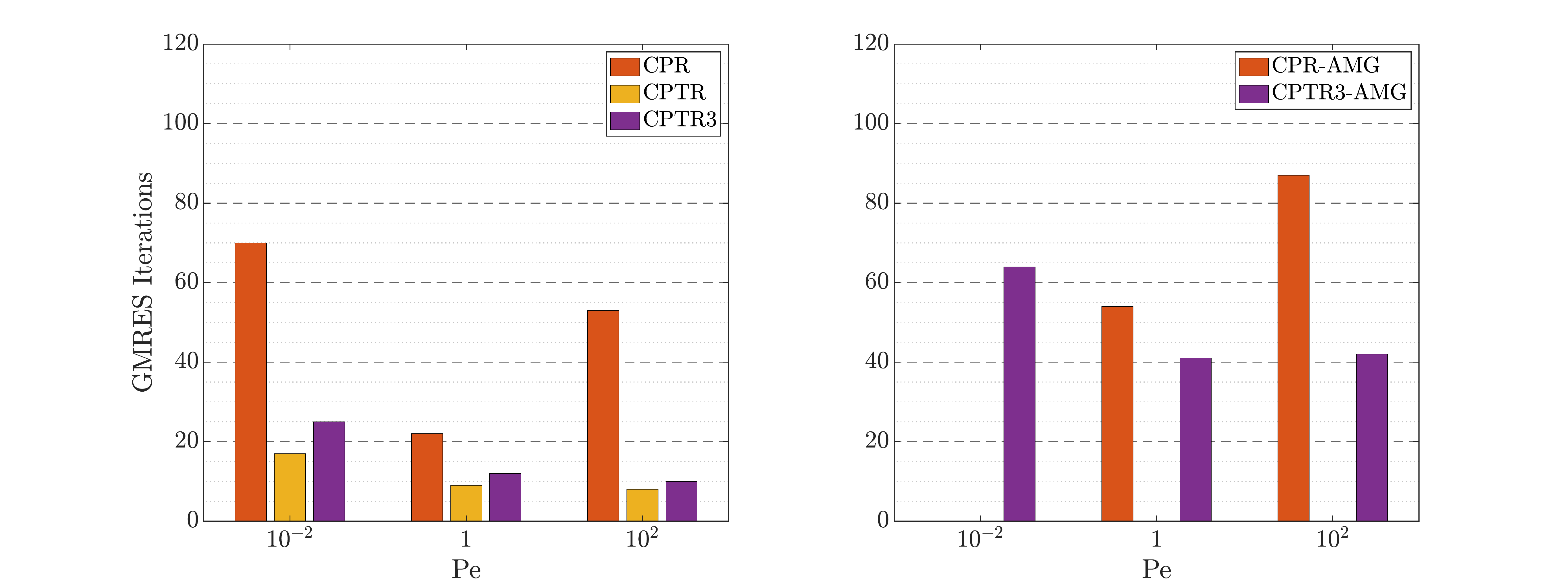}
    \caption{SPE10 bottom layer, results at $t=15$ min, using direct sub-solvers (left) and AMG sub-preconditioners (right).}
    \label{fig::MultPePaperBottom}
\end{figure}

We also test the bottom layer, representing a fluvial, channelized environment. The property variations are sharp between channels and levees, and the permeability spans more than six orders of magnitude. Figure \ref{fig::MultPePaperBottom} shows the performance of both new preconditioners using direct solvers for the subsystems (left) and using AMG  preconditioners (right). SPE10's bottom layer is one of the more challenging benchmark cases available regarding permeability variations. The coupling is even stronger than for the top layer, so the direct solver cases perform admirably.  The  subsystems become significantly more challenging, however, for AMG. Nevertheless, for the high P\'eclet case, CPR-AMG does not converge in under 200 iterations (it takes 227 iterations) but CPTR3-AMG makes that case tractable again (62 iterations). The average number of iterations across P\'eclet numbers for CPR-AMG is 123, while CPTR3-AMG achieves an average of only 49 iterations with a maximum of 62.

\subsection{In-Situ Combustion Homogeneous Case}
\label{subsec::isc}

One recovery method for heavy and extra-heavy oil is In-Situ Combustion (ISC), which involves chemical reactions (see \citet{Crookston79} and \citet{Coats80b} for details). Hot air injection is only efficient if one can leverage the chemical properties of the oil and oxidize part of it. Oxidation reactions are highly exothermic and allow the remainder of the oil to be displaced by lowering the viscosity. For reference, the oil we used throughout this study is 6.5$^\circ$API and the viscosity at reservoir conditions (48$^\circ$C) is close to 10,000 cP, rendering it virtually immobile for secondary processes like waterflooding. We consider the reaction scheme,
\begin{align*}
\label{eq::scheme}
\begin{split}
 & (1)\ \textrm{Oil} \rightarrow \alpha\textrm{Coke} \,, \\
 & (2)\ \textrm{Coke} + \beta\textrm{O$_2$} \rightarrow \gamma\textrm{CO$_2$} \,.
\end{split}
\end{align*}
Although simple, this scheme allows us to study the ideal combustion case of laying out the solid fuel (Coke) from the heaviest oil fraction (Oil) and burning the fuel in the presence of oxygen. We do not crack the oil and generate lighter products here, but simply use the enthalpy of reaction to lower the viscosity. In this case, one can completely consume both the oxygen and the fuel at the combustion front; the heat generated by combustion is located at the oxygen concentration front.

There are three main differences between the previous cases and this homogeneous combustion case. First, we generate heat inside the domain. The P\'eclet number is no longer sufficient to fully describe the thermal regime, but should be used in conjunction with the Damk\"ohler numbers. These compare reaction and advection ($\textrm{Da}_\textrm{I}$) and reaction and diffusion ($\textrm{Da}_{\textrm{II}}$). We will not perform a full dimensionless analysis in this case, but will present results for a typical ISC case using properties of unconsolidated sand for the thermal conductivity and parameters adapted from \citet{Dechelette06} for the reactions. If we still use the flow rate to compute the P\'eclet number, its value is $\textrm{Pe} = \mathcal{O}\left(10^{-1}\right)$, corresponding to a slightly diffusion dominated case. The second difference is that the reaction parameters (pre-exponential factor and activation energy) tend to be grid-sensitive and should be different across grid sizes to give similar temperature plateaus and front speeds \citep{Kovscek13,Nissen15}. Lastly, the beginning of the simulation is the ignition step, which is an unusual regime that we will not study in detail. To ensure the comparability of results, we ensure all cases for all grid sizes fall in the developed combustion regime and that the temperature profiles are similar.


\begin{figure}[t]
    \centering
    \includegraphics[width=\textwidth]{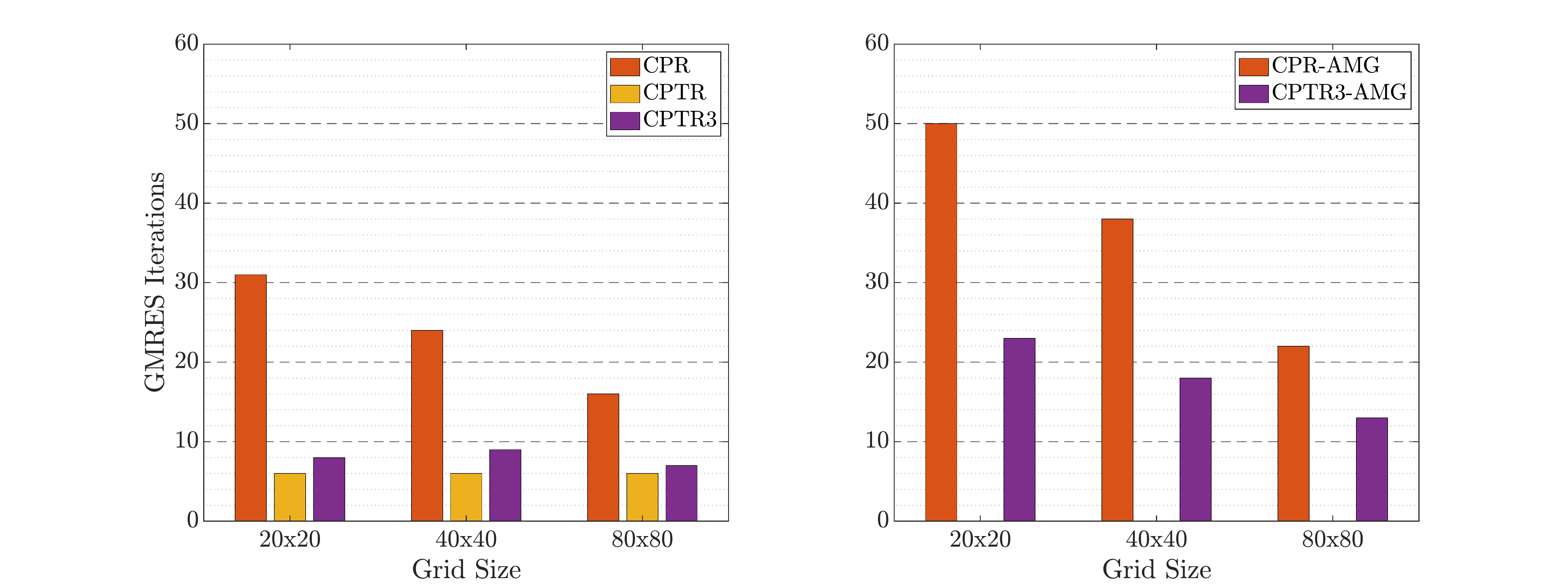}
    \caption{GMRES iterations results for the ISC case using a direct solver (left) an AMG preconditioner for all subsystems (right).}
    \label{fig::resultsISC25}
\end{figure}


Figure \ref{fig::resultsISC25} shows the results using  direct solvers (left) and the AMG preconditioners (right). For ISC, the maximum time step of 10 seconds cannot be sustained by the non-linear solver. We observe a lower number of GMRES iterations with finer grids, but caution should be taken in comparing results across mesh refinement levels.  When the grid size decreases, a smaller time step is required for nonlinear convergence.  Therefore, comparisons should only be made at a fixed mesh resolution, without inferring trends across resolutions. The CPR solver struggles, partly because of significant thermal diffusion. Another important effect is the much stronger coupling between the temperature and the saturations/mole fractions due to the reaction terms. The combustion front consumes oil and oxygen to release carbon dioxide, along with a large amount of energy. Failing to capture that coupling leads to significantly higher number of iterations for CPR. CPTR and CPTR3 outperform CPR by 81\% and 74\%, respectively. We observe a significant increase in the number of iterations for all grid sizes when using AMG preconditioners. Both the temperature, and especially the pressure subsystems, are much more challenging in the ISC case. Even with those more challenging subproblems, both CPTR-AMG and CPTR3-AMG still outperform CPR-AMG by around 55\% for the $40\times40$ and $80\times80$ grid sizes. In the context of thermal-compositional-reactive simulations, the non-linear solver is likely to struggle with the very tight coupling and the multi-scale nature of the physical phenomena. With many time-step cuts and small time steps, the number of linear solves per non-linear iterations will grow significantly, making it all the more important and valuable to be able to rely on a robust and fast preconditioning technique. CPTR3-AMG achieves that goal on our ISC test case, keeping the number of iterations below 22.

\section{Conclusion and Discussion}
\label{sec::concl}

In this paper we have presented two new multi-stage preconditioners for thermal-compositional(-reactive) flow in porous media. Using a dedicated treatment of temperature, either in a pseudo-two-stage (CPTR) or three-stage (CPTR3) fashion, the convergence of GMRES was significantly improved for all of our tests cases. The reduction in the number of iterations is 40-85\% compared to the industry standard CPR two-stage method, and at least one order of magnitude compared to a single-stage ILU(0) preconditioner on practical grid sizes for both homogeneous and heterogeneous cases.

The sensitivity of the proposed methods to the thermal regime, described by the thermal P\'eclet number in the absence of reactions, is greatly reduced and yields a much more robust preconditioner for a variety of reservoir simulation conditions. The use of AMG for the pressure and temperature subsystems and an ILU(0) final stage in CPTR3 shows good performance across all cases.

An interesting direction for future work is to apply a system AMG strategy for the elliptic subsystem of CPTR. Promising results have been obtained using SAMG for flow and transport problems \citep{Gries15}, and BoomerAMG for linear elasticity problems \cite{Baker11}. \citet{Roy19b} showed very promising results for dead-oil water injection using BoomerAMG. We have also only explored moderate grid sizes and serial solvers, but scalable, parallel implementations of these methods is an important next step. This work has also tackled just one element of improving the computational speed and robustness of thermal-compositional-reactive simulations, which exhibit many well-known computational bottlenecks. 

\section*{Acknowledgements}
\label{sec::acknow}

The authors want to thank Prof. Hamdi A. Tchelepi and Prof. Margot G. Gerritsen for interesting discussions and suggestions. Funding was provided by TOTAL S.A. through the FC-MAELSTROM project. Portions of this work were performed under the auspices of the U.S. Department of
Energy by Lawrence Livermore National Laboratory under Contract DE-AC52-07-NA27344.



\bibliographystyle{elsarticle-num-names}

\bibliography{bib_Nicola}



\end{document}